\documentclass[12pt]{amsart}
\usepackage{amsmath}
\usepackage{amssymb}
\usepackage{stmaryrd}
\usepackage{epsfig,color,esint}
\usepackage{blindtext}
\usepackage{enumerate}
\usepackage{enumitem}
\usepackage{hyperref}
\usepackage{graphicx}
\usepackage{url}
\usepackage{bbm}
\usepackage{nicefrac,mathtools}
\usepackage{mathrsfs}

\usepackage{comment}
\usepackage{pgfplots}
\usepackage{tikz}
\usetikzlibrary{arrows,shapes,trees,backgrounds}

\theoremstyle{definition}

\newtheorem{theorem}{Theorem}[section]

\newtheorem{proposition}[theorem]{Proposition}
\newtheorem{lemma}[theorem]{Lemma}
\newtheorem{remark}[theorem]{Remark}
\newtheorem{corollary}[theorem]{Corollary}

\newtheorem*{acknowledgements}{Acknowledgements}

\newtheorem*{remark*}{Remark}

\numberwithin{equation}{section}

\newcommand{\lb}[1]{\langle#1\rangle}

\newcommand{\mc}{\mathcal}
\newcommand{\mb}{\mathbb}

\newcommand{\eps}{\varepsilon}
\newcommand{\wti}{\widetilde}

\def\p{\partial}

\newcommand{\R}{\mathbb{R}}
\newcommand{\J}{\mathcal{J}}

\def\S{\mathbb S}
\def\K{\mathcal K}

\newcommand{\M}{\mathcal{M}}

\newcommand{\N}{\mathcal{N}}

\newcommand{\E}{\mathcal{E}}

\newcommand{\W}{\mathcal{W}}

\newcommand{\dist}{\operatorname{dist}}

\def\Om{\Omega}
\def\ol{\overline}
\def\D{\nabla}
\def\r{\vec{\hskip1pt r}}

\DeclareMathOperator{\vol}{Vol}

\def\R{\mathbb R}
\def\S{\mathbb S}
\def\J{\mathcal J}
\def\M{\mathcal M}
\def\N{\mathcal N}

\def\K{\mathcal K}
\def\P{\mathcal P}

\def\Cl{\text{Cl}}

\def\Om{\Omega}
\def\ol{\overline}
\def\D{\nabla}

\def\s{\sigma_{\S^n}}
\def\Int{\text{Int}\,}

\def\e{{\bf e}}
\def\ms{\mathcal S}

\newcommand{\beq}{\begin{equation}}
\newcommand{\eeq}{\end{equation}}
\newcommand{\beqs}{\begin{eqnarray*}}
\newcommand{\eeqs}{\end{eqnarray*}}
\newcommand{\beqn}{\begin{eqnarray}}
\newcommand{\eeqn}{\end{eqnarray}}
\newcommand{\beqa}{\begin{array}}
\newcommand{\eeqa}{\end{array}}

\def\lan{\langle}
\def\ran{\rangle}

\makeatletter
\@namedef{subjclassname@2020}{\textup{2020} Mathematics Subject Classification}
\makeatother

\title[$L_p$-Minkowski problem]{The $L_p$-Minkowski problem with  \\  super-critical exponents}

\author{Qiang Guang}
\address{Mathematical Sciences Institute, The Australian
         National University, Canberra, ACT 2601, Australia.}
\email{qiang.guang@anu.edu.au}

\author{Qi-Rui Li}
\address{School of Mathematical Sciences, Zhejiang University, Hangzhou 310027, China}
\email{qi-rui.li@zju.edu.cn}

\author{Xu-Jia Wang}
\address{Mathematical Sciences Institute, The Australian
         National University, Canberra, ACT 2601, Australia.}
\email{xu-jia.wang@anu.edu.au}


\begin{document}

\begin{abstract}
The $L_p$-Minkowski problem deals with the existence of closed convex hypersurfaces in 
$\R^{n+1}$ with prescribed $p$-area measures.
It extends the classical Minkowski problem and embraces several important geometric and physical applications.
Existence of solutions has been obtained in the sub-critical case $p>-n-1$, 
but the problem remains widely open in the super-critical case $p<-n-1$.
In this paper we introduce new ideas to solve the problem for all the super-critical exponents.
A crucial ingredient in our proof is a topological method 
based on the calculation of the homology of a topological space of ellipsoids. 

\end{abstract}


\maketitle
\baselineskip16.8pt
\parskip3.5pt

\section{Introduction}

A central problem in convex geometry is the characterisation of  geometric measures 
for convex bodies in the Euclidean space  $\R^{n+1}$.
The best-known example is the classical Minkowski problem, 
which was a major impetus for the development of fully nonlinear PDEs.
In the last three decades,
a focus of research in convex geometry is the $L_p$-Minkowski problem introduced by Lutwak \cite{Lut93}.
It includes the classical Minkowski problem ($p=1$), 
the logarithmic Minkowski problem ($p=0$), 
and the centro-affine Minkowski problem and elliptic affine spheres ($p=-n-1$) as special cases  \cite{BLYZ13, CW06}.
The $L_p$-Minkowski problem was derived from the Brunn-Minkowski theory, and
research of the problem paved the way for further development of this theory \cite{CLYZ09, Kol20, LYZ02}.
The $L_p$-Minkowski problem also plays a significant role in other applications.
Of particular interest is that it describes self-similar solutions to Gauss curvature flows \cite{An99, An03, BCD, Chow}, 
and its projective invariance in the case $p=-n-1$ makes it fundamental in image processing \cite{AGLM, AST}.

Let $\K_o$ denote the set of closed convex bodies in $\R^{n+1}$ with the origin in the interior.
For any $\Om\in\K_o$ and $p\in \mb{R}$, its $p$-area measure is defined as $d\ms_p = u^{1-p}d\ms$ \cite{Lut93},
where $u$ is the support function of $\Omega$ and $\ms$ is the classical surface area measure of $\Om$.
Given a finite non-negative Borel measure $\mu$ on the unit sphere $\S^n$,
the $L_p$-Minkowski problem asks for the existence of solutions $\Om\in\K_o$ such that its $p$-area measure 
coincides with the given measure $\mu$. If $d\mu = f d\s$ for a density function $f$ on $\S^n$, 
then the $L_p$-Minkowski problem can be formulated as finding
solutions to the   Monge-Amp\`ere equation
\beq\label{equ:lp}
\det(\D^2u+uI) =f u^{p-1}  \ \ \text{on} \ \S^n,
\eeq
where $\D$ denotes the covariant derivative with respect to an orthonormal frame on $\S^n$, and $I$ is the identity matrix. 

The last three decades have witnessed a great progress in the study of the $L_p$-Minkowski problem.
The problem can be divided into three cases.

\begin{itemize}
\item   In the sub-critical case $p>-n-1$ (with respect to the Blaschke-Santal\'o inequality),
the existence of solutions was obtained in \cite {CW06}. 
However, there is no uniform estimate for equation  \eqref{equ:lp} when $p<0$ \cite{JLW15}, 
and there may exist infinitely many solutions  when $p<-n$ \cite{HLW16, Li19}.
When $p=0$, \eqref{equ:lp} is called the logarithm Minkowski problem;
necessary and sufficient conditions for the existence of solutions were obtained in \cite{BLYZ13}
when the prescribed measure is an even Borel measure.
For $p\ge 1$, the existence and regularity of solutions were obtained in \cite{CW06, LO95}.

\item In the critical case $p=-n-1$, equation \eqref{equ:lp} is called the centro-affine Minkowski problem \cite{CW06}.
The quantity $u^{n+2}\det (\D^2u+uI)$ is invariant under projective transforms 
and plays a key role in affine geometry.
For example when $f=1$,   \eqref{equ:lp} is the equation for affine elliptic spheres \cite{Sch14}.
The projective invariance also makes it of great interest to  image processing \cite{AGLM}.
A Kazdan-Warner type condition \cite{CW06} implies that \eqref{equ:lp} admits no solutions
for a general positive function $f$. 
There are many works dealing with the critical case 
\cite{ACW,  JLZ16, JWW, Li19,  Lu20}.
However, results on the existence and multiplicity of solutions in this case are far from being satisfactory.
The main difficulty is that the normalisation of blow-up sequences does not lead to a unique limit model.

\item
In the super-critical case $p< -n-1$, there are some results in the one dimensional case.
In \cite{DZ12, SL15} the authors obtained the existence of $\frac \pi k$-periodic ($k$-fold symmetry, $k\ge 2$) convex solutions.
In  \cite {An03}, Andrews proved that 
when $f=1$ and $p\in [-7,-2)$,  a convex solution to \eqref{equ:lp} must be a circle; 
when $p<-7$, a convex solution to \eqref{equ:lp} is either the circle, or
a curve with $k$-fold symmetry.
In high dimensions, Zhu \cite {Zhu17} proved the existence of solutions when $f$ is a discrete measure
with no essential subspaces,
but we are unaware of any existence results when $f$ is a function.

\end{itemize}

There are many related research works on the $L_p$-Minkowski problem  \cite{BBCY, CHZ, CLZ19, HLYZ16, ZX}.
It is interesting to compare equation \eqref{equ:lp} with the semi-linear elliptic equation
\beq\label{equ:sc}
 - \Delta_{g_0}u+ c_nR_{g_0} u = f(x) u^{\gamma} \  \ \text{on} \ M,
\eeq
where $(M, g_0)$ is an $n$-dimensional Riemannian manifold, $c_n$ is a constant depending only on $n$, 
and $R_{g_0}$ is the scalar curvature.
There is a vast body of literature on equation  \eqref{equ:sc}.
In the sub-critical case $1<\gamma< \frac {n+2}{n-2}$,
there is a uniform estimate for solutions to 
\eqref{equ:sc}, and one can obtain the existence of non-trivial solutions under suitable conditions. 
In the critical case $\gamma= \frac {n+2}{n-2}$,  \eqref{equ:sc} is the prescribing scalar curvature equation.
In particular, it is Nirenberg's problem when $n=2$ and $M=\S^2$, and the Yamabe problem when $f=1$.
In this case, there is a very rich phenomena on the existence and multiplicity of solutions, and
one can find many significant results  \cite{Bre08, BM09, KMS09, Sch}.
In the super-critical case $\gamma > \frac {n+2}{n-2}$, 
numerous attempts have been made for the existence of non-trivial solutions to  \eqref{equ:sc} but the solution
was obtained only in some special cases.

Comparing with  \eqref{equ:sc}, we find that equation  \eqref{equ:lp} is more complicated.
There is no uniform estimate for  \eqref{equ:lp} in the sub-critical case.
There is no solutions in general in the critical case by the Kazdan-Warner type condition, 
and much less is known about sufficient conditions for the existence of solutions.
Therefore, 
one would not expect a complete resolution for the existence of solutions to \eqref{equ:lp} in the super-critical case.
Surprisingly,  we find that the $L_p$-Minkowski problem \eqref{equ:lp} admits a solution for all $p$ in the super-critical range, 
without any additional conditions, and thus completely resolve the existence problem.

\begin{theorem}\label{thm:main}
Suppose that $p<-n-1$.
Let $f$ be a positive and $C^{1,1}$-smooth function on $\S^{n}$.
Then there is a uniformly convex, $C^{3,\alpha}$-smooth and positive solution to \eqref{equ:lp},
where $\alpha\in (0,1)$.
\end{theorem}

By approximation, we also obtain the existence of solutions when $f$ is a non-smooth function. 

\begin{corollary}\label{thm:mainb}
Suppose that $p<-n-1$ and $f$ is a function on $\S^n$ such that $1/c_0\le f\le c_0$ for some constant $c_0>1$.
Then there is a strictly convex, $C^{1,\alpha}$-smooth and positive weak solution to \eqref{equ:lp} for some  $\alpha\in(0,1)$.
\end{corollary}

We point out that the condition $f>0$ in Theorem \ref{thm:main} and
Corollary \ref{thm:mainb} cannot be relaxed to $f\ge 0$.
Indeed, there exist functions $f$ which are positive  except at the north and south poles,
such that equation \eqref{equ:lp} admits no solutions  \cite [Theorem 1.4]{Du21}.

It is well-known that the Monge-Amp\`ere equation is of divergence form and 
equation \eqref{equ:lp} is the Euler equation of the following functional  (for $p\ne0$) for convex bodies $\Omega\in \K_o$
\beq\label{functional}
\J(\Om) = \vol(\Om) -\frac1p \int_{\S^n}f u^p d\s .
\eeq
Therefore, a natural approach to the $L_p$-Minkowski problem is 
to combine the variational method with the Gauss curvature flow \cite {BIS19, CW00, Iva,LSW20}.
In this paper we will employ the following Gauss curvature flow:
\begin{equation}\label{flow}
\frac{\p X}{\p t } (x,t)=-f(\nu)K(x,t)\lb{X,\nu}^p \nu+X(x,t),
\end{equation}
where $X(\cdot, t)$ is a parametrisation of the evolving convex hypersurfaces $\M_t$, 
$\nu$ and $K$ are respectively the unit outward normal and Gauss curvature of $\M_t$. 
We will show that the functional \eqref{functional} in non-increasing under the flow \eqref{flow} (Lemma \ref{lem:J}).

The main difficulty  is the lack of uniform estimate for the problem.
The uniform estimate is the key estimate for many geometric problems such as the Yamabe problem 
\cite{Sch} or Calabi‘s conjecture \cite {Yau}.
The $L_p$-Minkowski problem has been extensively studied in the past three decades, and
various techniques have been developed to establish the uniform estimate for the $L_p$-Minkowski problem
and the associated Gauss curvature flow, but none of them applies to the super-critical case. 

To overcome the difficulty, our strategy is to use a topological method to find a special initial condition 
such that the evolving hypersurfaces $\M_t=\p\Om_t$ satisfies 
\beq \label{BrR}
B_r(0)\subset \Om_t\subset B_R(0),
\eeq
for positive constants $R\ge r>0$ indenpendent of $t$,
where $B_r(x)$ denotes a closed ball of radius $r$ centred at $x$.  
Once the solution satisfies such a $C^0$-estimates,  
one can establish the second derivative estimates, and higher regularity follows from Krylov's regularity theory.
Hence by the monotonicity  of the functional \eqref{functional},  the flow converges to a solution of \eqref{equ:lp}.

Therefore, the key point in the argument is to find the special initial hypersurface. A crucial ingredient in achieving this goal is to compute the homology for a class of ellipsoids centred at the origin.
Let us outline the main ideas of the proof  below.

For any  convex body $\Om$ in $\R^{n+1}$, 
it is well known that there is a unique ellipsoid $E(\Om)$, called John's minimum ellipsoid \cite{Sch14},
which achieves the minimal volume among all ellipsoids containing $\Om$, such that
$${\Small\text{$\frac{1}{n+1}$}}E(\Om)\subset \Om\subset E(\Om). $$
Let $r_1(\Om)\le r_2(\Om)\le\cdots\le r_{n+1}(\Om)$ be the lengths of semi-axes of $E(\Om)$.
Denote $e_\M = e_\Om=\frac{r_{n+1}(\Om)}{r_1(\Om)}$
the eccentricity of $\M:=\p \Omega$ (or the eccentricity of $\Om$). 
We will first prove the following property. 
\begin{itemize}
\item [(P):]
For any given constant $A >\J(B_1(0))$, 
if one of the quantities $e_\Om$, $\vol(\Om)$, $[\vol(\Om)]^{-1}$, and  $[\dist(O,\p\Om)]^{-1}$ is sufficiently large,
we have $\J(\Omega)\geq A$
(Lemmas \ref{lem:J:d}-\ref{lem:J:e}).
\end{itemize}

Denote by $\mathcal A_I$ the set of ellipsoids $E$ such that 
the origin $O\in E$, $e_E\in[1,\bar{e}]$,  and
${{\bar v}}\le \vol(E)\le1/{{\bar v}} $, where  $\bar e$ is a large constant and ${{\bar v}}$ is a small constant.
$\mathcal A_I$ is a metric space under the Hausdorff distance.
For any ellipsoid $E\in \mathcal A_I$, let $\M_E(t)$ be the solution to the flow \eqref{flow} with initial condition $E$.
By the above property (P), $\M_{E}(t)$ has uniformly bounded eccentricity and volume,
and $\dist(O,\M_{E}(t))$ is uniformly bounded from zero if
\beq \label{JMA}
\J(\M_E(t))\le A\ \ \ \forall \ t \ge 0.
\eeq
Now, our focus is to prove that at any given time $t_0>0$, 
there exists an initial $E_0 \in \mathcal A_I$ such that the mimimum ellipsoid of $\M_{E_0}(t_0)$
is the unit ball centred at the origin (Lemma \ref{lem:ball}), 
 thus validating the condition \eqref{BrR} (as a result of Lemma \ref{lem:good}).

If to the contrary there is no such an initial $E_0$, 
we will construct a continuous map $T: \mc{A}_I \to \mathcal P$ which is the identity map on $\mc{P}$, where $\P$ is the boundary of $\mc{A}_I$ in the topological space of all ellipsoids.
This implies the existence of an injection from the homology group of $\mathcal P$ to that of $\mathcal A_I$.
As a consequence, $\mathcal P$ has trivial homology since $\mathcal A_I$ is contractible (Lemma \ref{lem:contr:A}).
By Proposition \ref{thm:topA} this leads to $H_k(\E\times\S^n) = H_k(\E)\oplus H_k(\S^n)$,
where $\E$ is introduced in \eqref{EEE}.
We thus reach a contradiction by the K\"unneth formula and Theorem \ref{thm:top} if we take $k =\frac{n(n+1)}{2}+2n-1$. 
This topological fixed-point argument is the main novelty in this paper.
A crucial ingredient in the argument is the computation of the homology of the class $\E$ of ellipsoids.

To complete the proof, we choose a sequence $t_k\to \infty$ and let $E_k\in \mathcal A_I$ be the initial condition 
such that the minimum ellipsoid of $\M_{E_k}(t_k)$ is  the unit ball.
By the Blaschke selection theorem, $E_k$ sub-converges to $E_*\in \mathcal A_I$.
It follows by the above property (P) that the Gauss curvature flow \eqref{flow} with  initial condition $E_*$
satisfies \eqref{BrR}.
Hence,  the flow \eqref{flow} starting from $E_*$ converges to a solution of \eqref{equ:lp} as $t\to\infty$.

The paper is organised as follows.
In Section \ref{sec:2}, we derive some a priori estimates for 
the functional \eqref{functional} and the Gauss curvature flow \eqref{flow}.
In Section \ref{sec:3}, we prove the main results (Theorem \ref{thm:main} and Corollary \ref{thm:mainb}),
assuming Proposition \ref{thm:topA}, Theorem \ref{thm:top}, and Theorem \ref{thm:estimates} temporarily. 
The proofs of Proposition \ref{thm:topA} and Theorem \ref{thm:top} will be given in Sections \ref{sec:top}, and the proof of Theorem  \ref{thm:estimates} will be given in Section \ref{sec:ap}.

The topological method introduced in this paper enables us to find an initial condition such that the solution has uniform estimate for all time $t$. The uniform estimate is the most difficult part for many geometric and analysis problems. 
This method can be adapted to other geometric problems, such as the $L_p$ dual Minkowski problem \cite{HLYZ16, LYZ18},
and the dual centro-affine Minkowski problem \cite {JLW18}
and more general prescribing curvature problems. 
We will study these problems separately.



\section{A priori estimates}\label{sec:2}

Let $\M$ be a smooth, closed, and uniformly convex hypersurface in $\R^{n+1}$.
The support function of $\M$ is given by
\beqs
u(x) = \lan x,\nu^{-1}_\M(x) \ran, \ \ \forall \ x\in\S^n,
\eeqs
where $\nu_\M:\M\to \S^n$ is the Gauss map and $\nu^{-1}_\M$ is its inverse,
namely, $\nu_\M^{-1}(x)$ is the point $z(x)\in \M$ such that the unit outer normal of $\M$ at $z(x)$ is equal to $x$.
It is well known that $\nu^{-1}_\M(x) = u(x)x+\D u(x)$ 
and the Gauss curvature of $\M$ at $\nu_{\M}^{-1}(x)$ is given by
\beq\label{s2 e1}
K = 1/\det (u_{ij}+u\delta_{ij}),
\eeq
where $u_{ij}:=\D^2_{ij}u$. 
This implies that the $p$-area measure of $\M$ is given by \cite{Lut93}
\beqs
d\mathcal S_p = u^{1-p}\det (\D^2u+uI)d\s .
\eeqs
Hence, the $L_p$-Minkowski problem is equivalent to solving  equation \eqref{equ:lp}.

Denote by $\Cl(\M)$ the convex body enclosed by $\M$.
When no confusion arises we may abuse the notation $\M$ for $\Cl(\M)$,
such as writing the functional $\J(\text{Cl}(\M))$  as  $\J(\M)$.
Assume that $\Cl(\M)\in\mathcal K_o$. Let $r$ be the radial function of $\M$, which is given by
\beq\label{rad funct}
r(\xi) = \max\{\lambda : \lambda\xi\in \Cl(\M)\} \ \ \forall \ \xi\in\S^n.
\eeq
Then 
\beq\label{equ:vol}
\vol(\Cl(\M)) = \frac{1}{n+1}\int_{\S^n} r^{n+1}d\s.
\eeq
Denote $\r(\xi) = r(\xi)\xi$. 
We also define the radial Gauss mapping by  
\beqs
\mathscr A_{\M}(\xi) = \nu_\M(\r(\xi)) \ \  \forall \ \xi\in\S^n. 
\eeqs

Let  $\M_t$ be a solution to the flow \eqref{flow} and $X(\cdot,t)$ be its parametrisation.
Consider the new parametrisation
\beqs
\ol X(x,t) = X(\nu_{\M_t}^{-1}(x),t).
\eeqs
It is straightforward to compute
\beqs
\frac{\p\ol X}{\p t} = \sum_i\frac{\p X}{\p z^i}\frac{\p (\nu_{\M_t}^{-1})_i}{\p t} + \frac{\p X}{\p t}.
\eeqs
Since the first term on the right hand side is tangential, 
taking inner product with the unit outer normal of $\M_t$ gives that
\beqs
\p_tu(x,t) = \big\lan x,\p_t \ol X(x,t) \big\ran = \big\lan x,\p_t X(\nu_{\M_t}^{-1}(x),t)  \big\ran.
\eeqs
Hence by \eqref{s2 e1}, the flow \eqref{flow} can be expressed as 
\begin{equation}\label{flow:u}
\p_t u(x,t)=-\frac{f(x) u^p(x)}{\det(\D^2u+uI)}+u(x,t).
\end{equation}

We next show the monotonicity of the functional \eqref{functional} under the flow \eqref{flow}.

\begin{lemma}\label{lem:J}
Suppose $\M_t$, $t\in [0,T)$, is a solution to the flow \eqref{flow} in $\K_o$. Then 
\[\frac{d}{dt}\J(\Omega_t)\geq 0,\]
where $\Om_t=\Cl(\M_t)$. 
Moreover, the equality holds if and only if $\M_t$ satisfies \eqref{equ:lp}.
\end{lemma}

\begin{proof}
The following formulas can be found in \cite{LSW20}:
\beq\label{lem:J e1}
{\begin{split}
  {\Small{\text{$ \frac{\p_t r}{r}(\xi,t) $}}} & ={\Small{\text{$ \frac{\p_t u}{u}$}}} (\mathscr A_{\M_t}(\xi),t) ,\\
   |\text{Jac}\mathscr A|(\xi)  & =  {\Small{\text{$\frac{r^{n+1}K(\r(\xi,t))}{u(\mathscr A_{\M_t}(\xi))} $}}} ,
  \end{split}}
\eeq
where $\text{Jac}\mathscr A$ is the Jacobian of the radial Gauss mapping.

By virtue of \eqref{s2 e1}-\eqref{lem:J e1}, we obtain
\beqs
\frac{d}{dt}\J(\Omega_t)
&=& - \int_{\S^n} fu^{p-1}\p_tu(x)d\s(x) + \int_{\S^n} r^{n}\p_t r(\xi)d\s(\xi)\\
&=& \int_{\S^n}\big(  {\Small{\text{$  \frac{1}{K} $}}} - fu^{p-1} \big)\p_tu(x)d\s(x) \\
&=&\int_{\S^n}\big( {\Small{\text{$ \frac1K $}}} -fu^{p-1}\big)^2  uKd\s \ge 0.
\eeqs
Clearly, the equality 
$\frac{d}{dt}\J(\Omega_t) =0$ holds if and only if $u(\cdot,t)$ satisfies \eqref{equ:lp}.
\end{proof}

The proof of Lemma \ref{lem:J} also verifies that  \eqref{equ:lp} is the Euler-Lagrangian equation of the 
functional \eqref{functional}.

 \bigskip 
 
\subsection{Properties of the functional \eqref{functional}}
{~}

Next, we prove the  property (P) stated in the introduction.

\begin{lemma}\label{lem:J:d}
Suppose that $p<-n-1$ and $1/c_0\le f\le c_0$ for some $c_0\geq 1$. 
For any given constant $A>0$, there exists a small constant ${d_0}>0$ depending only on $n$, $p$, $c_0$ and $A$ 
such that if $\Omega\in \K_o$  satisfies 
$\text{dist}(O,\p \Omega)\in (0, {d_0})$, then
$\J(\Omega)>A.$
\end{lemma}

\begin{proof}
Denote by $d=\text{dist}(O,\p \Omega)>0$. 
Take $x_0\in \mb{S}^n$ such that 
\[u(x_0)=\min_{\mb{S}^n} u=d,\]
where $u$ is the support function of $\Omega$. 
Let $E$ be the minimum ellipsoid of $\Omega$. 
We choose the coordinates such that 
\[E-\zeta_E=\Big\{z\in \mb{R}^{n+1}:  {\Small\text{$ \sum_{i=1}^{n+1} $}} \, \frac{z_i^2}{a_i^2}\leq 1\Big\},\]
and
\[x_0\cdot \e_{n+1}=\max\{|x_0\cdot \e_i|: 1\leq i\leq n+1\} , \]
where $\zeta_E$ is the center of $E$. 
This implies that $x_0\cdot \e_{n+1}\geq c_n$.
We use $c_n$ to denote a constant which depends only on $n$, 
but it may change from line to line.  

Let $w(x)=u(x)+u(-x)$, $x\in \mb{S}^n$, be the width function of $\Omega$. 
Since the ball $B_d(0)$ is contained in $\Omega$ and $\frac{1}{n+1}E\subset \Omega\subset E$, we have
\[d\leq \min_{\mb{S}^n} w\leq c_n a_{n+1} \quad\text{and }\quad w(\e_i)\leq c_n a_i. \]
This yields that 
\begin{equation}\label{e:J:d1}
\J(\Omega)>\vol(\Omega)
      \geq c_n  {\Small\text{$ \prod_{i=1}^{n+1} $}} a_i
     \geq c_n d  {\Small\text{$ \prod_{i=1}^n$}}  w(\e_i) 
     \geq c_n d  {\Small\text{$ \prod_{i=1}^n$}} u(\e_i).
\end{equation}
 
Next, we consider the set $\Omega_{n+1}^*=\Omega^*\cap L$, 
where $\Omega^*$ is the polar dual of $\Omega$ and $L=\{z\in \mb{R}^{n+1}: z\cdot \e_{n+1}=0\}$. 
Let $r^*$ be the radial function of $\Omega^*$. 
Since the origin $O$ and points $r^*(\e_i)\e_i$, $i=1, \cdots, n$,  are contained in $\Omega_{n+1}^*$, 
their convex hull is an $n$-dimensional convex set in $\Omega_{n+1}^*$,
namely,
\beqs
\mathcal C=: \text{convex hull of }\{O, r^*(\e_1)\e_1, \cdots, r^*(\e_n)\e_n\} \subset \Omega_{n+1}^*.
\eeqs
Let $V$ be the cone in $\mb{R}^{n+1}$ with base $\mathcal C$ and vertex $p_0=r^*(x_0)x_0$.
Since $r^*(x)=\frac{1}{u(x)}$ and $u(x_0)=d$, the height of the cone $V$ (in the direction of $\e_{n+1}$) satisfies 
\begin{equation}\label{e:J:d3}
r^*(x_0)x_0\cdot \e_{n+1}=\frac{x_0\cdot \e_{n+1}}{d}\geq \frac{c_n}{d}.
\end{equation}
Consider the following subset of $V$
\[V'=\Big\{z\in V: z_{n+1}\geq \frac{r^*(x_0)}{2} x_0\cdot \e_{n+1}\Big \}.\]
By \eqref{e:J:d3}, we have 
\begin{equation}\label{e:J:d4}
|z|\geq \frac{c_n}{2 d} \quad \forall\ z\in V'.
\end{equation}
In view of $V'\subset \Omega^*$, $f\geq 1/c_0$ and \eqref{e:J:d4}, we have  
\beqs
{\begin{split}
\J(\Omega) & \geq -\frac{1}{p}\int_{\mb{S}^n} u^p f\geq \frac{1}{c_0}\int_{\Omega^*} |z|^{-p-n-1}dz \\
                   & \geq \frac{1}{c_0} \int_{V'} |z|^{-p-n-1}dz\geq \frac{c_n\vol(V')}{c_0d^{-p-n-1}}. 
                   \end{split}}
\eeqs
Since
\beqs
{\begin{split}
   & \vol(V')\geq \frac{c_n}{d} \vol(\mathcal C),\\
   & \vol(\mathcal C)\geq c_n  {\Small\text{$ \prod_{i=1}^n$}} r^*(\e_i)
   =c_n \Big[  {\Small\text{$ \prod_{i=1}^n $}} u(\e_i)\Big]^{-1},
                   \end{split}}
\eeqs
we further obtain that 
\begin{equation}\label{e:J:d5}
\J(\Omega)\geq \frac{c_n}{c_0d^{-p-n}}\Big[  {\Small\text{$ \prod_{i=1}^n$}} u(\e_i)\Big]^{-1}. 
\end{equation}
Combining \eqref{e:J:d1} and \eqref{e:J:d5}, we have 
\[\big[\J(\Omega)\big]^2\geq \frac{c_n}{c_0d^{-p-n-1}}.\]
Since $-p-n-1>0$, we see that  $\J(\Omega)>A$ if  $d$ is sufficiently small. 
\end{proof}

\begin{lemma}\label{lem:J:vol} 
Suppose that $p<-n-1$ and $1/c_0\le f\le c_0$ for some $c_0\geq 1$.
For any given constant $A>1$, there exists a  small constant $v>0$ depending only on $n, p,c_0$, and $A$,  
such that if $\Omega\in \K_o$  satisfies
either  $\vol(\Omega)\leq v$ or $ \vol(\Omega)\geq v^{-1}$,
then
$
\J(\Omega)> A.
$
\end{lemma}
\begin{proof}
If $\vol(\Omega)\geq v^{-1}$, by definition we have $\J(\Om) > \vol(\Omega)\geq v^{-1}>A$ by taking $v$ small. 

If $\vol(\Omega)\leq v$. Denote $d=\text{dist}(O,\p \Omega)$. 
Since the ball $B_{d}(0)\subset \Omega$, we have 
\[ {v}\geq \vol(\Omega)\geq \vol(B_d)=\vol(B_1)d^{n+1}.\]
Hence, if $v$ is sufficiently small, then $d<{d_0}$, where ${d_0}>0$ is the constant given by Lemma \ref{lem:J:d}.
Therefore, we have  $\J(\Omega)> A$ by Lemma \ref{lem:J:d}.
\end{proof}

\begin{lemma}\label{lem:J:e}
Suppose that $p<-n-1$ and $1/c_0\le f\le c_0$ for some $c_0\geq 1$.
For any given constant $A>0$, 
there exists a large constant $e>1$ depending only on $n,p,c_0$ and $A$, such that
if $\Omega\in \K_o$ satisfies  $e_\Om\ge e$, we have
$
\J(\Omega)> A.
$
\end{lemma}

\begin{proof}
For the given constant $A$, let  ${d_0}$  be the constant determined by Lemma \ref{lem:J:d}.  
We assume that $d=\text{dist}(O,\p \Omega)\geq {d_0}$; otherwise, we are done by Lemma \ref{lem:J:d}.

Let $E$ be the minimum ellipsoid of $\Omega$ with semi-axes $a_1\leq \cdots \leq a_{n+1}$. Note that $\frac{a_{n+1}}{a_1}=e_\Omega$.
Since $B_{d_0}\subset \Omega$ and $\Omega\subset E$, we obtain that 
\begin{equation}\label{e:J:e1}
{d_0}\leq d\leq a_1.
\end{equation}
Noting also that $\frac{1}{n+1}E\subset \Omega$, we have 
\[\J(\Omega)> \vol(\Omega)\geq c_n \vol(E)=c_n {\Small\text{$ \prod_{i=1}^{n+1} $}} a_i\geq c_n e_\Omega {d_0}^{n+1}.\]
Clearly, we have $\J(\Omega)>A$ if $e_\Omega$ is sufficiently large. 
\end{proof}

\vskip10pt

\subsection{A priori estimates for the parabolic equation \eqref{flow:u}}
{~}

In this subsection, we state the a priori estimates for the solution $u$,
assuming the uniform estimate for $u$.

\begin{theorem}\label{thm:estimates}
Let $f$ be a positive and $C^{1,1}$-smooth function on $\S^n$.
Let $u(\cdot,t)$ be a positive, smooth and uniformly convex solution to \eqref{flow:u}, $t\in[0,T)$. 
Assume that 
\beq\label{thm:estimates:C0}
{\begin{split}
   & 1/C_0\le u(x,t)\le C_0,\\
   & \ \ \ |\D u|(x,t)\le C_0,
  \end{split}} 
\eeq
for all $(x,t)\in\S^n\times[0,T)$. Then
\beq\label{D2u}
C^{-1}I\le  (\D^2u+uI)(x,t)\le  CI \ \ \forall \ (x,t)\in\S^n\times[0,T),
\eeq
where $C$ is a positive constant depending only on $n, p, C_0$, $\min_{\S^n}f$, $\| f \|_{C^{1,1}(\S^n)}$,
and the initial condition $u(\cdot,0)$.
\end{theorem}

The proof of Theorem \ref{thm:estimates} is based on proper choice of auxiliary functions, and will be given in 
Section \ref{sec:ap}.

By the second derivative estimates  \eqref{D2u},
equation \eqref{flow:u} becomes uniformly parabolic.
Hence, by Krylov's regularity theory,  we have the $C^{3,\alpha}$ estimates for the solution $u$. Namely
\beq \label {D3u}
  \|u(\cdot,t) \|_{C^{3,\alpha}(\S^{n})} \le C  \ \ \forall \ (x,t)\in\S^n\times[0,T),
\eeq 
for any given $\alpha\in (0, 1)$, where the constant $C$ depends only on 
$ \alpha, n,p,\min_{\S^n}f, \|f\|_{C^{1,1}(\S^n)}$, and the initial condition $u(\cdot,0)$. 
By the a priori estimates \eqref{D3u}, 
we have the longtime existence of solutions to the flow \eqref{flow}, provided that $u$ satisfies \eqref{thm:estimates:C0}.

\begin{theorem}\label{thm:longtime}
Let $f$ be a positive and $C^{1,1}$-smooth function on $\S^n$.
Let $T_{\max}$ be the maximal time such that 
$u(\cdot, t)$ is a positive, $C^{3,\alpha}$-smooth, and uniformly convex solution to \eqref{flow:u} on $[0,T_{\max})$.
If \eqref{thm:estimates:C0} holds for all the time $t\in[0,T_{\max})$, 
then $T_{\max}=\infty$ and $u$ satisfies the estimates
\eqref{D2u} and \eqref{D3u}. 
\end{theorem}

\begin{remark} \label{r2.7}
Let $\M(t)_{|\, t\in [0,T_{\max})}$ be a solution to \eqref{flow}.
By Lemmas \ref{lem:J:d}, \ref{lem:J:vol} and \ref{lem:J:e}, if $\J(\M(t))<A$ for some constant $A$ independent of $t$,
then there exists positive constants $e, v, d$ depending on $A$, but independent of $t$, such that
\beq\label{evd}
e_{\M(t)}\le e, \ \ \ v\le \vol(\Om(t))\le v^{-1},\ \text{ and }\ B_{d}(0)\subset \Om(t),
\eeq
where $\Om(t)$ is the convex body enclosed by $\M(t)$. 

From \eqref{evd} one infers that \eqref{thm:estimates:C0} holds.
Hence  the a priori estimates \eqref{D2u} and \eqref{D3u} hold, 
and one has the long-time existence of solution (Theorem \ref{thm:longtime}).
Therefore, for the a priori estimates \eqref{D2u}, \eqref{D3u} and the the long-time existence of solution,
all we need is that the condition $\J(\M(t))<A$ holds for some constant $A$.
\end{remark}

\section{Proof of Theorem \ref{thm:main} }\label{sec:3}

In this section, we show how to select an initial hypersurface $\N_0$, such that 
the flow \eqref{flow} deforms $\N_0$ to  a solution of \eqref{equ:lp}.

The initial hypersurface $\N_0$ is an ellipsoid and  will be chosen  by a topological method.
In the proof of the existence of $\N_0$, 
a key step is the computation of  the homology groups of a special class of ellipsoids.
The homology groups will be given in Proposition \ref{thm:topA} and Theorem \ref{thm:top}, 
whose proofs are postponed to the next section.

Denote
\beq\label{s3 A0}
A_0= 2\Big(-\frac{\|f\|_{L^1(\S^n)}}{p[2(n+1)]^p}+2^{n+1}\vol(B_1)\Big).
\eeq
where $B_1=B_1(0)$ is the unit ball in $\R^{n+1}$ centred at the origin.  
Recall that $p<-n-1$.  Hence for any
$\Om\in\mathcal K_o$ with $\frac{1}{2(n+1)}B_1 \subseteq\Om\subseteq 2B_1$, we have
\beqs
\J(\Om)\le \frac12A_0 . 
\ \ \eeqs
In particular, if the minimum ellipsoid of $\Om$ is $B_1$, then $\frac{1}{n+1}B_1\subset\Om\subset B_1$ and hence
\beqn\label{s3.1 e1}
\J(\Om) \le -\frac{1}{p(n+1)^p}\int_{\S^n}f+\vol(B_1)\le  \frac12A_0.
\eeqn

\bigskip
\subsection{A modified flow of \eqref{flow}}\label{subsec:flow}
{~}
 
We introduce a modified flow of \eqref{flow} such that for any initial condition, 
the solution exists for all  time $t\ge 0$.
The purpose of introducing this modified flow is for the convenience of  later discussion.

For a closed, smooth and uniformly convex hypersurface $\N$ such that $\Om_0 = \Cl(\N)\in \K_o$,
we define a family of time-depending hypersurfaces $\bar{\M}_\N(t)$ with initial condition $\N$ as follows:
 
\begin{itemize}

\item [{\Small\text{$\bullet$}}]  If $\mc{J}(\M_\N(t)) < A_0$ for all time $t\ge 0$, let 
$\bar{\M}_\N(t)= \M_\N(t)$ for all $t\ge 0$, where 
$\M_\N(t)$ is the solution to \eqref{flow}.
 We point out that, by Remark \ref{r2.7}, the solution $\M_\N(t)$ exists as long as $\mc{J}(\M_\N(t)) $ is finite.

\item [{\Small\text{$\bullet$}}]  If $\J(\N)<A_0$, and $\mc{J}(\M_\N(t))$ reaches $A_0$ at the first time $t_0>0$, 
we define
\[
\bar{\M}_\N(t)=\left\{ 
\begin{split}
& \M_\N(t), \ & \ \text{if}\ 0\leq t<t_0,\\
& \M_\N(t_0),\ & \text{if}\ \   t\geq t_0.
\end{split}
\right.
\]

\vskip3pt

\item [{\Small\text{$\bullet$}}]  If $\J(\N)\geq A_0$, we let $\bar{\M}_\N(t)\equiv \N$ for all $t\geq 0$. That is, the solution is stationary.
\end{itemize}

For convenience, we call $\bar{\M}_\N(t)$ a modified flow of \eqref{flow}.

\begin{remark} Apparently $\Cl(\bar{\M}_{\N}(t))\in\K_o$ and it is easy to verify the following properties.
\begin{itemize}
\item[$\diamond$] $\bar{\M}_\N(t)$ is defined for all time $t\ge 0$, and by Lemma \ref{lem:J},
$\J(\bar{\M}_\N(t))$ is  non-decreasing.

\item[$\diamond$] By Lemma \ref{lem:J:d}, if $\dist(O,\N)$ is very small, then $\bar{\M}_\N(t)\equiv \N$ $\forall \ t\ge 0$.

\item[$\diamond$]By Lemma \ref{lem:J:vol}, 
if $\vol(\Omega_0)$ is sufficiently large or small, 
then  $\bar{\M}_\N(t)\equiv \N$ $\forall \ t\ge 0$.

\item[$\diamond$]By Lemma  \ref{lem:J:e}, 
if $e_{\Om_0}$ is sufficiently large, then  $\bar{\M}_\N(t)\equiv \N$ $\forall \ t\ge 0$.

\item[$\diamond$] We have $\J(\bar{\M}_\N(t)) \le \max \{A_0, \J(\N)\}$ $\forall\ t\ge 0$. 

\item[$\diamond$] By the a priori estimates, 
$\bar{\M}_\N(t)$ is smooth for any fixed time $t$, and Lipschitz continuous in time $t$. 

\end{itemize}
\end{remark}

\subsection{A special class of ellipsoids $\mathcal A_I$}

\begin{lemma}\label{lem:e-v}
For the constant $A_0$ given by \eqref{s3 A0}, 
there exist sufficiently small constants $\bar v$ and $\bar d$,
and a sufficiently large constant  $\bar e$,
such that for any $\Om\in \mathcal K_o$,
\begin{itemize}
\item[(i)] if $\text{dist}(O,\p \Omega)\leq \bar{d}$, then $\J(\Om)> A_0$;  
\item[(ii)]  if $e_\Om\ge \bar e$, then $\J(\Om)> A_0$;
\item[(iii)]  if $\vol(\Om)\le \bar v$ or $\vol(\Om)\ge [(n+1)^{n+1}\bar v]^{-1}$, then $\J(\Om)>A_0$.
\end{itemize}
\end{lemma}

\begin{proof}
This is an immediate consequence of Lemmas  \ref{lem:J:d}, \ref{lem:J:vol} and \ref{lem:J:e}. See also Remark \ref{r2.7}.
\end{proof}

Fix the constants $\bar d, \bar v$ and $ \bar e$ as in Lemma \ref{lem:e-v}.
We introduce the following notations:
\begin{itemize}
\item [{\Small\text{$\bullet$}}] $\K$ is the collection of all non-empty, compact and convex sets in $\R^{n+1}$
equipped with the Hausdorff distance,  such that $\K$ is a metric space. 

\item [{\Small\text{$\bullet$}}] $\bar \K_o$ is the closure of $\K_o$ in $\K$.

\item [{\Small\text{$\bullet$}}]  $\K_e$ is the subset of $\K_o$ which consists of all origin-symmetric convex bodies.

\item [{\Small\text{$\bullet$}}]  Denote 
\beqs
&&\mathcal A_I = \big\{ E\in \bar \K_o \ \text{is an ellipsoid in } \mb{R}^{n+1}: 
        \ {{\bar v}}\le \vol(E)\le {{1/\bar v}} \ \text{and}\ e_E\le \bar e\big\},\\
&&\hat{\mathcal A}=\{E\in\mathcal A_I: \ \vol(E)=\omega_n \ \text{and} \ e_E\in[1,\bar e]\}, \quad\text{where }\omega_n=\vol(B_1),\\
&& {{\mathcal A}}
=\{E\in\hat{\mathcal A}: \ \text{either} \ e_E=\bar e\ \text{or} \ \dist(O,\p E)=0\}.
\eeqs

\item [{\Small\text{$\bullet$}}]  
To calculate the homology of  $\mathcal A_I$, we also introduce
\beq \label{EEE}
{\begin{split}
 \E_I & = \mathcal A_I\cap \K_e ,\\  
 \hat\E & = \hat{\mathcal A}\cap \K_e ,\\
  \E &= \mathcal A\cap \K_e.
  \end{split}}
\eeq

\end{itemize}

\noindent 
The sets $\mathcal A_I, \hat{\mathcal A}, {{\mathcal A}}$ and $ \E_I, \hat\E, \E$ are all closed. In particular, $\E$ (resp. $\mc{A}$) is the boundary of $\hat{\E}$ (resp. $\hat{\mc{A}}$) in the space of all ellipsoids in $\K_e$ (resp. $\K$) with unit ball volume.



For any $E\in\mathcal E$, let $\mathcal R$ be a rotation such that
\beqs
\mathcal R(E)=\Big\{ (x_1,\cdots,x_{n+1})\in\R^{n+1}: \   {\Small\text{$\sum_{i=1}^{n+1}$}}\, \frac{x_i^2}{a_i^2}\le 1\Big\}.
\eeqs
Since $\vol(E) = \omega_n$, we have $\prod_{i=1}^{n+1}a_i = 1$. For $s\in[0,1]$, denote
\beqs
b_i(s) = (1-s)+sa_i\ \ \text{and}  \ \ \hat a_i(s)=\frac{b_i(s)}{[\prod_{i=1}^{n+1} b_i(s)]^{\frac{1}{n+1}}}.
\eeqs
We then obtain an ellipsoid $\hat E(s)$ such that
\beq \label{hatE}
\mathcal R(\hat E(s)) = \Big\{(x_1,\cdots,x_{n+1})\in\R^{n+1}: 
    \  {\Small\text{$ \sum_{i=1}^{n+1}$}}\, \frac{x_i^2}{ (\hat a_i(s))^2} \le 1\Big\}.
\eeq
In particular, $\hat E(0) =  B_1$, $\hat E(1) = E $. 
The set  $\{\hat E(s)\}_{s\in[0,1]}$ is a path in $\K_e$ connecting $ B_1$ to $E$, 
and satisfies $\vol(\hat E(s)) = \omega_n$ for all $s\in[0,1]$.
As a result,
\beq\label{def hat E}
\hat{\mathcal E} = \{ \hat E(s): \ E\in \mathcal E , s\in[0,1]\},
\eeq
namely,  $\hat{\mathcal E}$ is the union of all such paths, and
\beq\label{def EI}
\mathcal E_I = \big\{ [\tau/\omega_n]^{\frac{1}{n+1}} \hat E: \ \hat E\in\hat{\mathcal E}, \ \tau\in [{{\bar v}},{{1/\bar v}}]\big \}.
\eeq

\begin{lemma}\label{lem:contr}
Both $\mathcal E_I$ and $\hat{\mathcal E}$ are contractible. 
Hence the homology 
\beq\label{HHe}
H_k(\mathcal E_I) = H_k(\hat{\mathcal E}) = 0\ \ \ \forall \ k\ge1 .
\eeq
\end{lemma}

\begin{proof} 
To see that $\hat{\mathcal E}$ is contractible, by \eqref{def hat E} it suffices to notice that 
$\big(\hat E(s),t\big)\to \hat E((1-t)s)$, where $t\in [0, 1]$, is a deformation retraction from
$\hat{\mathcal E}\to \{B_1\}$.
Hence $\hat{\mathcal E}$ are contractible. 

Similarly, by \eqref{def EI},
\[\Big(\Big[\frac{\tau}{\omega_n}\Big]^{\frac{1}{n+1}} \hat E ,t\Big)\mapsto 
\Big[\frac{(1-t)\tau+t\omega_n}{\omega_n}\Big]^{\frac{1}{n+1}} \hat E \]
is a deformation retraction from $\mathcal E_I$ to $\hat{\mathcal E} $. 
As $\hat{\mathcal E} $ is contractible, one sees that $\mathcal E_I$ is also contractible. 
\end{proof}

We can combine the two deformation retractions in the above proof
and obtain a new one from $\mathcal E_I$ to $ \{B_1\}$ as follows
\beq\label{comb}
\eta:\ \ \Big(\Big[\frac{\tau}{\omega_n}\Big]^{\frac{1}{n+1}} \hat E(s) ,t\Big)\mapsto 
\Big[\frac{(1-t)\tau+t\omega_n}{\omega_n}\Big]^{\frac{1}{n+1}} \hat E((1-t)s) ,
\eeq
where $s\in [0, 1]$, $\tau\in [{{\bar v}}, {{1/\bar v}}]$. 

In the following, 
the notation $\mathcal H\simeq\mathcal H'$ means two metric spaces $\mathcal H, \mathcal H'$ are homeomorphic.

\begin{lemma}\label{lem:contr:A}
We have $\mathcal A_I\simeq \mathcal E_I\times B_1$,
and $\hat{\mathcal A}\simeq \hat\E\times B_1$ (with the product topology). It follows that both $\mathcal A_I$ and $\hat{\mathcal A}$ are contractible,
and  hence the homology 
\beq\label{HHa}
H_k(\mathcal A_I) = H_k(\hat{\mathcal A}) = 0\ \ \ \forall \ k\ge1. 
\eeq
\end{lemma}
\begin{proof}
For any given ellipsoid $E\in \K_e$, the map
\beqs
\vec{\hskip1pt\varrho}_E({{y}})=
\left\{ 
\begin{split}
 \frac{{{y}}}{r_E({{y}} / |{{y}} | )}, \quad  & \text{if} \ {{y}}\in E\setminus\{0\},\\
 0, \quad\quad\quad &\text{if} \ {{y}}=0 , 
\end{split}
\right.
\eeqs
defines a homeomorphism between $E$ and $B_1$,
where $r_E$ is the radial function of $E$ (see \eqref{rad funct}).
Denote by $\vec{\hskip1pt\varrho}_E^*$ the inverse of $\vec{\hskip1pt\varrho}_E$.

We now show $\mathcal A_I\simeq\E_I\times B_1$.
For $E\in \mc{A}_I$, let $\zeta_E$ be its centre. 
Then $E_o=: E-\zeta_E\in \E_I$. 
We define a map $\phi: \mc{A}_I\to \E_I\times B_1$ by
\beq\label{map phi}
\phi(E)=\big(E-\zeta_E, \vec{\hskip1pt\varrho}_{E_o}(\zeta_E)\big).
\eeq
Its inverse $\phi^*: \mc{E}_I\times B_1\to \mc{A}_I$ is given by 
\beq\label{map phi*}
\phi^*(E, {{y}})=E+\vec{\hskip1pt\varrho}^*_{E}({{y}}).
\eeq
It is easy to verify that $\phi^*\circ \phi=id_{\mathcal A_I}$ , 
$\phi \circ \phi^*=id_{\E_I\times B_1}$ and both $\phi$ and $\phi^*$ are continuous.

Restricting $\phi$ and $\phi^*$ to $\hat{\mathcal A}$ and $\hat\E\times B_1$ respectively,
we see that $\hat{\mathcal A}\simeq \hat\E\times B_1$. \end{proof}

Denote
\beq\label{def:P}
\mathcal P
=\big\{E\in\mathcal A_I: \text{either}\ \vol(E)={{\bar v}},\ \text{or} \ \vol(E)={{1/\bar v}}, \ \text{or} \ e_E=\bar e, \ \text{or} \ O\in \p E\big\}.
\eeq
Note that $\mc{P}$ is the boundary of  $\mathcal A_I$ if we regard $\mathcal A_I$ as a set in the topological space of all ellipsoids.

\begin{lemma}\label{lem:Psi}
There is a retraction $\Psi$ from $\mathcal A_I\setminus \{B_1\}$ to $\mathcal P$. 
Namely, $\Psi: \mathcal A_I\setminus \{B_1\}\to \mathcal P$ is continuous and $\Psi|_{ \mathcal P}=id$.
\end{lemma}

\begin{proof}
The retraction $\Psi$ can be constructed as follows. 

\noindent
{1.}
By the map $\phi$ in \eqref{map phi}, we have
$\mathcal A_I\setminus\{B_1\}\simeq \E_I\times B_1\setminus (\{B_1\}\times\{0\})$.

\vskip5pt
\noindent
{2.}
Let $\p \E_I$ be the topological boundary of $\E_I$, regarding $\E_I$ as a set in the space of all ellipsoids in $\K_e$.
Then $\p \E_I=\{E\in\E_I:\text{either}\ \vol(E) = {{\bar v}}, \ \text{or} \ \vol(E)={{1/\bar v}}, \ \text{or} \ e_E=\bar e\}$.
In this step we define a continuous map 
$\psi: (\E_I\times B_1)\setminus( \{B_1\}\times\{0\})\to (\p\E_I\times B_1 ) \cup( \E_I \times \p B_1)$.

Recall the deformation retraction $\eta$ from $\mathcal E_I$ to $ \{B_1\}$ in \eqref{comb}.
For any $t \ne 1$, $\eta(\cdot, t)$ defines a homeomorphism between $\E_I$ and $\eta(\E_I, t)$.
For any $E \in \E_I\setminus \{B_1\} $, let 
$$t_E=\sup \{ t\in [0, 1]:\ \text{there exists $E' \in \E_I$ such that $\eta(E', t)=E$} \}. $$
We have $t_E<1$.
Since $\E_I$ is closed, there exists $\tilde E\in \p \E_I\subset\E_I$ such that $\eta(\tilde E, t_E)=E$. 
For any given $t\in [0, 1)$, $\eta$ also defines a homeomorphism between
$\p E_I$ and $\eta(\p E_I, t)$, and $\eta(\cdot, 0)$ is the identity map on $\p E_I$.
Define a map $\psi_1:\ \E_I\setminus \{B_1\} \to \p \E_I$ by letting  $\psi_1(E)= \tilde E$, where $\tilde{E}$ satisfies $\eta(\tilde E, t_E)=E$.
Then $\psi_1$ is a retraction from $\E_I\setminus \{B_1\}$ to $\p \E_I$ such that $\psi_1=id$ on $\p \E_I$.

Let $\psi_2(x)=\frac {x}{|x|}$. 
Then $\psi_2$ is a retraction $B_1\setminus \{0\}$ to $\p B_1$ such that $\psi_2=id$ on $\p B_1$.

Combining $\psi_1$ and $\psi_2$ we can define the map $\psi$  by letting
$$
\psi:\  (E, x) \to \left\{
{\begin{split}
 & \big( \psi_1(E) ,  \frac{x}{1-t_E}\big)  \hskip90pt \text{if}\ 1-t_E\ge |x|,\\
 & \big(\eta\big(\psi_1(E), 1-\frac{1-t_E}{|x|}\big),  \psi_2(x)\big) \ \hskip15pt \text{if}\ 1-t_E<|x|.
\end{split}}
\right. $$

\noindent
{3.}
By the map $\phi^*$ in \eqref{map phi*}, we have $\p(\E_I\times B_1 ) \simeq \mathcal P$.

\vskip5pt
\noindent
{4.} Let $\Psi = \phi^*\circ\psi\circ\phi :\mathcal A_I\setminus \{B_1\}\to \mathcal P$.
Then $\Psi$ is a retraction from $\mathcal A_I\setminus \{B_1\}$ to $\mathcal P$.
\end{proof}

The following two results are crucial for our later argument, whose 
proofs are postponed to Section \ref{sec:top}.

\begin{proposition}\label{thm:topA}
We have the following results.
\begin{itemize}
\item[(i)] $H_{k+1}(\mathcal P) =H_k({{\mathcal A}})$ for all $k\ge1$.
\item[(ii)] There is a long exact sequence
\beqs
\cdots\to H_{k+1}({{\mathcal A}})\to H_k(\E\times\S^n)\to H_k(\E) \oplus H_k(\S^n) \to H_k({{\mathcal A}})\to\cdots.
\eeqs
\end{itemize}
\end{proposition}

\begin{theorem}\label{thm:top}
Let $n^* = \frac{n(n+1)}{2}$.
The homology group $H_{n^*+n-1}(\mathcal E)=\mathbb Z$.
\end{theorem}

\bigskip

\subsection{Selection of a good initial condition}
{~}

As mentioned in the introduction, 
we use a topological method to prove the existence of a special initial condition 
such that the solution to the Gauss curvature flow \eqref{flow} satisfies the uniform estimate \eqref{BrR}.

We will employ the  modified  flow with initial data in $\mathcal A_I$.
For any ellipsoid $\N$ such that $\Cl(\N)\in \mathcal A_I$,
let $\bar{\M}_\N(t)$ be the solution to the modified flow,
with the constant $A_0$ given in \eqref{s3 A0}. 
We have the following properties
\begin{enumerate}
\item If $\Cl(\N)$ is close to $\mathcal P$ in Hausdorff distance or on $\mathcal P$, 
we have $\mc{J}(\N)\ge A_0$ and so $\bar{\M}_\N(t)\equiv \N$ for all $t$ (see Lemma \ref{lem:e-v}). 
\item If $\Cl(\N)$ is close to $B_1(0)$ in Hausdorff distance,  then $\J(\N)<A_0$.
\item By our definition of the modified flow,  if $\J (\bar{\M}_\N(t)) < A_0$ for all $t\geq 0$, then by Remark \ref{r2.7}, we have
\beq\label{s3.3 e0}
e_{\bar{\M}_\N(t)}\le \bar{e}, \ \ \ {{\bar v}}\le \vol(\bar\Om_{\N}(t))\le {{1/\bar v}},\ \text{ and }\ 
B_{\bar d}(0)\subset \bar\Om_{\N}(t)\ \forall\ t\ge 0,
\eeq
where  $\bar\Om_{\N}(t) = \Cl(\bar{\M}_\N(t))$,  the convex body enclosed by $\bar{\M}_\N(t)$.
Here the bar over $\Omega$ means that $\bar\Om_{\N}(t)$ is the convex body enclosed 
by the modified flow $\bar{\M}_\N(t)$,  not the closure of $\Om_{\N}(t)$. 
In our notation, $\Om_{\N}$ is a closed convex body.
\end{enumerate}
With  these properties, 
we can prove the following key lemma.
 
\begin{lemma}\label{lem:ball}
For every $t>0$, there exists $\N=\N_t$ with $\Cl(\N)\in \mc{A}_I$, 
such that the minimum ellipsoid of $\bar{\M}_{\N}(t)$ is the unit ball $ B_1(0)$.
\end{lemma}

\begin{proof}
Suppose by contradiction that there exists $t'>0$ such that, for any $\Om\in\mathcal A_I$,
$E_\N(t')\ne  B_1$, where $\N=\p\Om$ and $E_\N(t')$ is the minimum ellipsoid of $\Om_{\N}(t'):=\Cl(\bar{\M}_\N(t'))$. We have  
\beq\label{lem:ball e1}
E_\N(t')\in\mathcal A_I.
\eeq
This is obvious when $\Omega_{\N}(t)\equiv \Omega$ for all $t$. If $\Omega_{\N}(t)$ is not identical to $\Omega$, then \eqref{lem:ball e1} follows from $\frac{1}{n+1}E_\N(t')\subset \Om_{\N}(t')\subset E_\N(t')$, \eqref{s3.3 e0} and Lemma \ref{lem:e-v}.


Hence, we can define a continuous map $T: \mathcal A_I\to \mathcal P$ by 
\beqs
\Om\in \mathcal A_I\mapsto E_\N(t')\in \mathcal A_I\setminus  \{B_1\}\mapsto \Psi(E_\N(t'))\in \mathcal P,
\eeqs
where $\Psi$ is given in Lemma \ref{lem:Psi}.
Note that when  $\Omega\in \mathcal P$, we have $\J(\Omega)\ge A_0$ and thus 
$E_\N(t')=E_\N(0)=\Omega$.  This implies  that $T|_{\mathcal P}=id_{\mathcal P}$.
Hence  $T$ is a retraction from $\mathcal A_I$ to $\mathcal P$,
and so there is an injection from $H_*(\mathcal P)$ to $H_*(\mathcal A_I)$.
By \eqref{HHa} we then have 
\beqs
H_k( \mathcal P) = 0 \ \text{for all}\ k\ge1.
\eeqs
It follows from  Proposition \ref{thm:topA} (ii) that
\beqs
H_{k}(\mathcal E\times\S^n) =H_k(\E)\oplus H_k(\S^n) \ \forall \ k\ge1.
\eeqs
Computing the left-hand side by the K\"unneth formula and using the fact $H_k(\S^n) = \mathbb Z$ if $k=0$ or $k=n$, 
and $H_k(\S^n)=0$ otherwise, 
we further obtain
\beq\label{FH}
H_k(\E)\oplus  H_{k-n}(\E) = H_k(\E)\oplus H_k(\S^n)\ \ \forall\ k\ge 1.
\eeq
However, this contradicts Theorem \ref{thm:top} by taking $k = n^*+2n-1$ in the above.
\end{proof}

\begin{remark}\label{r3.9}
Lemma \ref{lem:ball} asserts that for any $t>0$,
there is an initial hypersurface
$\N=\N_t$ such that the minimum ellipsoid $E_\N(t)$ is the unit ball $B_1(0)$,
even in the case when  the $L_p$-Minkowski problem \eqref{equ:lp} has a unique solution $\p B_R(z)$ for $R\ne 1$ and $z\ne 0$.
In fact, our topological argument implies that for any ellipsoid $E\in \mc{A}_I\setminus \mc{P}$,
and any given time $t>0$, there is an initial hypersurface $\N=\N_t$ such that $E_\N(t)=E$.
\end{remark}

\begin{remark}\label{r3.10}
To derive a contradiction from \eqref{FH},
we only need to show that there exists one nontrivial homology group among $H_k(\E)$ for all $k\ge 1$.
Theorem \ref{thm:top} asserts that $H_{n^*+n-1}(\mathcal E)=\mathbb Z$, 
which suffices for the proof of the key Lemma \ref{lem:ball}.
By \eqref{large k hom} below, $k=n^*+n-1$ is the largest integer such that $H_k(\E)\ne0$. 
We will not compute other homology groups of $\E$ in this paper, as they are not needed in our proof.
\end{remark}

In the following we prove the convergence of the flow \eqref{flow} with a specially chosen initial condition.
Take a sequence $t_k\to\infty$ and let $\N_k = \N_{t_k}$ be the initial data from Lemma \ref{lem:ball}. 
By our choice of $A_0$ (see \eqref{s3 A0} and \eqref{s3.1 e1}), Lemma \ref{lem:ball} implies that 
\beq\label{s3.3 e1}
\mc{J}(\bar{\M}_{\N_k}(t_k))\le\frac12A_0.
\eeq
Hence, by the monotonicity of the functional $\J$, we have
\beqs
\bar{\M}_{\N_k}(t)=\M_{\N_k}(t) \ \ \forall \ t\le t_k.
\eeqs
Since $\Cl(\N_k)\in\mathcal A_I$ and $B_{\bar d}(0) \subset \Cl(\N_k)$, by Blaschke's selection theorem, 
there is a subsequence of $\N_k$ which converges in Hausdorff distance to a limit $\N_*$ such that $\Cl(\N_*)\in\mathcal A_I$
and $B_{\bar d}\subset \Cl(\N_*)$.

Next, we  show that the flow \eqref{flow} starting from $\N_*$ satisfying $\J(\M_{\N_*}(t))<A_0$ for all $t$.

\begin{lemma}\label{lem:good} 
For any $t\geq 0$, we have \[\mc{J}(\bar{\M}_{\N_*}(t))\le \frac{3}{4}A_0.\]
Hence 
$$\bar{\M}_{\N_*}(t)=\M_{\N_*}(t) \ \ \forall \ t>0.$$
\end{lemma}

\begin{proof}
For any given  $t>0$, since $\N_k\to \N_*$ and $t_k\to\infty$, 
when $k$ is sufficiently large such that $t_k>t$, we have
\beqs
\mc{J}(\bar{\M}_{\N_*}(t))-\mc{J}(\M_{\N_k}(t))\leq \frac14A_0.
\eeqs
By the monotonicity of the functional $\J$,   
\beqs
\mc{J}(\M_{\N_k}(t))\le \mc{J}(\M_{\N_k}(t_k)).
\eeqs
Combining above two inequalities with \eqref{s3.3 e1}, we obtain that 
\beqs
 \mc{J}(\bar{\M}_{\N_*}(t))&=&\mc{J}(\bar{\M}_{\N_*}(t))-\mc{J}(\M_{\N_k}(t)) +\mc{J}(\M_{\N_k}(t))\\ 
&\le&\mc{J}(\bar{\M}_{\N_*}(t))-\mc{J}(\M_{\N_k}(t))+\mc{J}(\M_{\N_k}(t_k))\\
& \leq&\frac14A_0 + \frac12A_0=\frac{3}{4}A_0.
\eeqs
This completes the proof.
\end{proof}

\vskip10pt

\subsection{Convergence of the flow and existence of solutions to \eqref{equ:lp}}
{~}

Let $\Om_{\N_*}(t) = \Cl(\M_{\N_*}(t))$ and $u(\cdot,t)$ be its support function.
By Lemma \ref{lem:good},  ${\M}_{\N_*}(t)$ satisfies \eqref{s3.3 e0}.
Hence by  \eqref{s3.3 e0} we infer that
\beqs
\bar d\leq u(x, t)\leq   C
\ \ \forall \ (x, t)\in\S^n\times[0,\infty).
\eeqs
where $C=(n+1)/(\bar v \omega_{n-1}\bar d^n) $.
By the convexity, 
\beqs
|\D u(x, t)| \le \max_{\S^n} u(\cdot,t)\leq C \ \ \forall \ (x,t)\in\S^n\times[0,\infty).
\eeqs
Namely condition \eqref{thm:estimates:C0} holds.
By Subsection 2.2, 
we obtain the existence of solutions to \eqref{equ:lp} as follows.

\vskip10pt

\begin{proof}[Proof of Theorem \ref{thm:main}]
Denote $\M(t)=\M_{\N_*}(t)$ 
and $\mc{J}(t)=\mc{J}(\M(t))$.
By Lemma \ref{lem:J} and Lemma \ref{lem:good},  
\[\J(t)< A_0\ \text{and }\ \J'(t)\geq 0\ \ \forall\ t\ge0.\]
Therefore,
\[\int_0^\infty \mc{J}'(t)dt \le \limsup_{T\to\infty} \J(T)-\J(0)  \leq A_0.\]
This implies that there exists a sequence $t_i\to \infty$ such that 
\beqs
\J'(t_i)=\int_{\S^n} \Big[\Big(\frac1K - fu^{p-1}\Big)^2 uK \Big]\Big|_{t=t_i} d\s\to 0.
\eeqs
Passing to a subsequence, we obtain by the a priori estimates \eqref{D3u} that
$u(\cdot,t_i)\to u_\infty$ in $C^{3,\alpha}(\S^{n})$-topology and $u_\infty$ satisfies \eqref{equ:lp}.
\end{proof}

Corollary \ref{thm:mainb} follows from Theorem \ref{thm:main}  by an approximation argument.
To prove Corollary \ref{thm:mainb}, 
 we first point out that all arguments in Subsections 2.1 and 3.1-3.3 depend on $\inf f$ and $\sup f$
but are independent of the smoothness of $f$. 
Therefore, the constants $\bar e, {{\bar v}},  \bar d$ in \eqref{s3.3 e0} are independent of the smoothness of $f$.

\begin{proof}[Proof of Corollary \ref{thm:mainb}]
Choose a sequence of functions $f_j\in C^\infty(\S^n)$ such that $\inf_{\S^n} f\le f_j\le \sup_{\S^n} f$ 
and $f_j\to f$ a.e. (such as the mollifications of $f$).
By our proof of Theorem \ref{thm:main},  
there exists an initial convex hypersurface $\N_j\in \mathcal A_I$ such that the solution $\M_j(t):=\M_{\N_j}(t)$ to \eqref{flow}
converges to a solution $\M_j$ of \eqref{equ:lp} with $f=f_j$, 
and  $\M_j(t)$ satisfies \eqref{s3.3 e0}, uniformly for all $j$ and $t$.
Hence $\M_j$ satisfies the estimates in \eqref{s3.3 e0}.
Passing to a subsequence, we may assume that $\M_j$ converges in Hausdorff distance to a limit $\M$.
Then $\M$ satisfies   \eqref{s3.3 e0}.
By the weak convergence of the Monge-Amp\`ere equation, $\M$ is a weak solution of \eqref{equ:lp}.
By the regularity theory of the Monge-Amp\`ere equation, 
$\M$ is strictly convex and $C^{1,\alpha}$ smooth for some $\alpha\in (0, 1)$.
\end{proof}



\section{Proofs of Proposition \ref{thm:topA} and Theorem \ref{thm:top}}\label{sec:top}

In this section, we prove Proposition \ref{thm:topA} and Theorem \ref{thm:top}. 

\subsection{Proof of Proposition \ref{thm:topA}} \label{sec:topA}{~}

\begin{proof}[Part (i)]
{~}
Note that for any  $\tau\in[\bar v,\bar v^{-1}]$, we have 
\beqs
&& \{E\in\mathcal A_I:\vol (E) = \tau\}\simeq \hat{\mathcal A}, \\
&& \{E\in\mathcal A_I:\vol(E) = \tau;  \ e_E=\bar e, \ \text{or} \ O\in\p E\}\simeq \mathcal A.
\eeqs
Hence, $\mathcal P$ consists of three components (up to homeomorphism)
\beqs
\mathcal A\times[\bar v,\bar v^{-1}], \ \hat{\mathcal A}\times\{\bar v\}, \ \text{and} \ \hat{\mathcal A}\times\{\bar v^{-1}\}.
\eeqs
Since $\mathcal A\times[\bar v, \bar v^{-1}]$ is topologically a cylinder with base $\mathcal A$,  and $\mc{A}$ can be viewed as the boundary of $\hat{\mc{A}}$, we see that $\mathcal P$ can be viewed as attaching two copies of $\hat{\mathcal A}$ along the boundary of 
$\mathcal A\times[\bar v,\bar v^{-1}]$.
As $\hat{\mathcal A}$ is contractable, we conclude that $\mathcal P$ is homotopy-equivalent to $S\mathcal A$
(the suspension of $\mathcal A$),
which is the quotient of $\mathcal A\times[\bar v,\bar v^{-1}]$ obtained by collapsing $\mathcal A\times\{\bar v\}$ to one point
and $\mathcal A\times\{\bar v^{-1} \}$ to another point.
Hence the two spaces have the same homology
$H_*(\mathcal P) = H_*(S\mathcal A)$.
It is known from \cite{Hat} that 
\beqs
H_{k+1}(S\mathcal A) = H_k(\mathcal A), \ \text{for all} \ k\ge1.
\eeqs
This completes the proof.
\end{proof}

\vskip10pt
\begin{proof}[Part (ii)]

Let $\mathcal B$ be the set of unit balls such that the origin lies on the boundary of the ball.
Consider the following subspaces of ${\mathcal A}$:
\beqs
{{\mathcal A}}_1={{\mathcal A}}\setminus\mathcal B,\ {{\mathcal A}}_2={{\mathcal A}}\setminus \E,\
\text{and} \ {{\mathcal A}}_3={{\mathcal A}}_1\cap{{\mathcal A}}_2.
\eeqs
We have the Mayer-Vietories sequence for the decomposition ${{\mathcal A}}={{\mathcal A}}_1\cup{{\mathcal A}}_2$:
\beq\label{s4a:2:e1}
\cdots\to H_{k+1}({{\mathcal A}})\to H_k({{\mathcal A}}_3)
\to H_k({{\mathcal A}}_1) \oplus H_k({{\mathcal A}}_2) \to H_k({{\mathcal A}})\to\cdots.
\eeq
Let $L =( [0,1]\times\{1\})\cup(\{1\}\times[0,1])\subset  \R^2$.  
Denote
\beqs
L_1 = \{(s,\rho)\in L: s>0\}, \ L_2 = \{(s,\rho)\in L: \rho>0\},  \ \text{and}\ L_3=L_1\cap L_2.
\eeqs
Let $G_1:L_1\times[0,1]\to L_1$ be a strong deformation retraction from $L_1$ onto the point $(1,0)$; 
$G_2:L_2\times[0,1]\to L_2$ be a strong deformation retraction from $L_2$ onto $(0,1)$;
and $G_3:L_3\times[0,1]\to L_3$ be a strong deformation retraction from $L_3$ onto $(1,1)$.
Denote by $G_{i,1}$ and $G_{i,2}$ the components of the map $G_i$ such that
$G_i(s,\rho,t)=(G_{i,1},G_{i,2})(s,\rho,t)$, where $1\leq i\leq 3$.

Given $E\in {{\mathcal A}}$, we take $(E',\xi,s,\rho)\in\E\times\S^n\times L$ such that
$\phi(E) =  \big(\hat{E'}(s),\rho\xi\big)$ with $\phi$ being the map \eqref{map phi}.
Define $\mathcal G_i:\mathcal A_i\times[0,1]\to\mathcal A_i$ by
\beqs
\mathcal G_i(E,t) = \phi^*(\hat E'(G_{i,1}(s,\rho,t)), G_{i,2}(s,\rho,t)\xi).
\eeqs
Then $\mathcal G_1,\mathcal G_2$ and $\mathcal G_3$ are deformation retractions
from ${{\mathcal A}}_1$, ${{\mathcal A}}_2$ and ${{\mathcal A}}_3$ onto 
$\E$, $\mathcal B$ and $\mathcal A'$ respectively, where 
\beqs
\mathcal A' = \{E\in\mathcal A: e_E=\bar e \ \text{and} \ O\in\p E\}.
\eeqs 
Therefore,
\beqs
 H_*({{\mathcal A}}_1) = H_*(\E), \   H_*({{\mathcal A}}_2) = H_*(\mathcal B),
  \ \text{and} \ H_*({{\mathcal A}}_3) = H_*(\mathcal A').
\eeqs
Inserting these identities into \eqref{s4a:2:e1}, we obtain the following long exact sequence
\beq\label{s4a:2:e2}
\cdots\to H_{k+1}({{\mathcal A}})\to H_k(\mathcal A')
\to H_k(\E) \oplus H_k(\mathcal B) \to H_k({{\mathcal A}})\to\cdots.
\eeq
On the other hand, the maps $\phi$ and $\phi^*$ (see \eqref{map phi} and \eqref{map phi*}) yield the homeomorphisms
\beq\label{s4a:2:e3}
\mathcal B\simeq\S^n \ \text{and} \ \mathcal A'\simeq \E\times\S^n.
\eeq
Our conclusion follows by combining \eqref{s4a:2:e2} and \eqref{s4a:2:e3}.
\end{proof}

\vskip10pt

The remaining of this section is devoted to the proof of Theorem \ref{thm:top}.

\subsection{Proof of Theorem \ref{thm:top}}
{~}

In the following, we always take $n^*=\frac{n(n+1)}{2}$, which is the dimension of the Lie group $SO(n+1)$.

For $n=1$, the proof of Theorem \ref{thm:top} is straightforward.
In this case, the lengths of semi-axes $r_1=r_1(E)$ and $r_2=r_2(E)$ of $E\in\mathcal E$ satisfy
$r_1r_2=1$ and $r_2 = \bar er_1$. Hence,
\beqs
r_1(E)= \frac{1}{\sqrt{ \bar e}} \ \ \text{and} \ \ r_2(E) = \sqrt{\bar e}, \ \ \forall E\in \E.
\eeqs
Therefore, each element $E\in\E$ is determined by the major axis of $E$.  
It follows that $\E$ is homeomorphic to $RP^1$. As a result,
\beqs
H_{1}(\E) = H_1(RP^1) = \mathbb Z.
\eeqs
This proves Theorem \ref{thm:top} when $n=1$, since $n^*=\frac{n(n+1)}2 = 1$.
 
 \bigskip

In the following we deal with the case $n\ge2$.
First we introduce some notations.
Let
\beqs 
{\begin{split}
H_1 &= \big\{(x_2,\cdots,x_{n}) \in\R^{n-1}: \ x_2\ge 1\big\} ,\\
H_i &= \big\{(x_2,\cdots,x_{n}) \in\R^{n-1}: \ x_{i+1}\ge x_i\big\} ,\ \ i=2, \cdots, n-1,\\
H_n &= \big\{(x_2,\cdots,x_{n}) \in\R^{n-1}: \ \bar e\ge x_n\big\} .
\end{split}}
\eeqs
Then $H_i$ are closed half spaces of $\R^{n-1}$, $i=1,2,\cdots,n$. Denote 
\beq\label{s4 e1}
\Delta_{n-1}  = \bigcap_{i=1}^n H_i.
\eeq 
If $n=2$,  then $\Delta_{n-1}  =\big\{1\le x_2\le \bar e\big\} \subset\R$ is an interval. If $n=3$,  then $\Delta_{n-1}  =\big\{(x_2,x_3) \in\R^2: \  1\le x_2\le x_3\le \bar e\big\}$
is a triangle. For $n\ge 3$, $\Delta_{n-1}$ is an $(n-1)$-simplex in $\mb{R}^{n-1}$.

Denote by
\beqs
F_i=\p H_i\cap \Delta_{n-1},
\eeqs
a face of $\Delta_{n-1}$, $i=1, \cdots, n$.
We also denote by  $\Delta_{n-1}^{(i)}$ the subset of $\Delta_{n-1}$,
obtained by removing the face $F_i$ from $\Delta_{n-1}$, namely
\beq\label{s4 e1b}
\Delta_{n-1}^{(i)}= \Delta_{n-1}\setminus F_i.
\eeq

There is a natural projection $ \pi:\E\to\Delta_{n-1}$, given by
\beqs
 E\mapsto \pi (E) =\big (\tilde r_2(E),\cdots,\tilde r_n(E) \big),
\eeqs
where $r_i(E)$ are lengths of semi-axes of $E$ satisfying $r_1(E)\le r_2(E)\le \cdots\le r_{n+1}(E)$ and
\beqs
\tilde r_i(E) = \frac{r_i(E)}{r_1(E)}, \ \ \ i=2, \cdots, n.
\eeqs
Note that $r_{n+1}(E)/r_1(E)=\bar e$ is a fixed constant for all $E\in\E$.

The mapping $\pi$ can be written as a composition of two mappings $\pi_1$ and $\pi_2$. Namely,
\beqs
\begin{split}
\pi_1:\E\to \mathbb L_{n+1}, \ \ \ \ \  \ \ &\text{given by} \ E\mapsto (r_1(E),\cdots,r_{n+1}(E)),\\
\pi_2:\mathbb L_{n+1}\to \Delta_{n-1},  \ \  &\text{given by} \ 
(x_1,\cdots,x_{n+1})\mapsto \big(\frac{x_2}{x_1},\cdots,\frac{x_{n}}{x_1}\big),
\end{split}
\eeqs
where
\beqs
\mathbb L_{n+1} = \Big\{(x_1,\cdots,x_{n+1})\in\R^{n+1}: 0<x_1\le \cdots \le x_{n+1}, \  {\Small\text{$ \prod_{i=1}^{n+1}$}} x_i=1,\ 
x_{n+1} = \bar{e} x_1 \Big\}.
\eeqs
It is readily seen that $\pi_2$ is a bijection and $\pi = \pi_2\circ\pi_1$ is surjective.

\begin{remark}\label{r4.1}
For any given ellipsoid $E\in\E$, there is a unique  positive definite, unimodular matrix $A$
such that $E=\{x\in\R^{n+1}:\ x'\cdot Ax=1\}$.  
Let $\lambda_1\ge \cdots\ge \lambda_{n+1}$ be the eigenvalues of $A$.
Then $\pi_1(E)=(\lambda_1^{-\frac{1}{2}}, \cdots, \lambda_{n+1}^{-\frac{1}{2}})$. 
Hence, for any point (vector) $ \r =(r_1, \cdots, r_{n+1})\in\R^{n+1}$ 
such that $0<r_1\le  \cdots \le r_{n+1}$ and $\prod_{i=1}^{n+1} r_i=1$, 
we can define $\pi_1^{-1}(\r)$ as the set of ellipsoids $E\in\E$ such that
the lengths of the semi-axes of $E$ are equal to $r_1, \cdots, r_{n+1}$, namely, 
the eigenvalues of the ellipsoid matrix $A$ are equal to $r_1^{-2}, \cdots, r_{n+1}^{-2}$.
For convenience, we say that a component $r_i$ of the vector $\r$ is single if $r_{i-1}<r_i<r_{i+1}$,
and that a component $r_i$ has multiplicity $k$ ($k\ge 2$) if $r_{i-1}< r_i=\cdots=r_{i+k-1}<r_{i+k}$.
\end{remark}

For $1\le i\le n$, consider the following subsets of $\E$:
\beq\label{s4 e2}
\E_i = \pi^{-1}(\Delta^{(i)}_{n-1}). 
\eeq
 By our notation $\Delta_{n-1}^{(i)}$ in \eqref{s4 e1b}, $\E_i$ is a subset of $\E$
and can be  written as
\beqs
\E_i= \big\{E\in\E:  r_i(E)\ne r_{i+1}(E)\big\}.
\eeqs
For any $1\leq j_1<j_2\cdots <j_l\leq n$, we denote
\beq\label{s4 e3}
\E_{j_1; j_2; \cdots; j_l}=\bigcup_{s=1}^l\E_{j_s}
\quad \text{and} \quad
\mathcal W_{j_1; j_2; \cdots; j_l}=\bigcap_{s=1}^l \E_{j_s}.
\eeq
For brevity we write 
$ \E_{j_1 j_2 \cdots j_l}= \E_{j_1; j_2; \cdots; j_l}$ and 
$ \mathcal W_{j_1 j_2 \cdots j_l}=
\mathcal W_{j_1; j_2; \cdots; j_l}$. 
We see that
\beqs
&&\E_{j_1j_2\cdots j_l} =\big \{ E\in\E: r_i(E)\ne r_{i+1}(E) \ \text{for some} \ i=j_1,\cdots,j_l\big\},\\
&&\W_{j_1j_2\cdots j_l} =\big \{ E\in\E: r_i(E)\ne r_{i+1}(E) \ \text{for all} \ i=j_1,\cdots,j_l\big\}.
\eeqs

For convenience of the reader, we first prove Theorem \ref{thm:top} for lower dimensions $ n=2$ and $n=3$,
and then for higher dimensions. One may also skip subsections \ref{subsec:n2} and \ref{subsec:n3}, 
and go through subsection \ref{subsec:n} for general case directly.
Our method is based on dividing $\E$ into suitable parts and employing the Mayer-Vietoris sequences \cite{Hat}.

\subsection{Dimension $n=2$}\label{subsec:n2}
{~}

For $n=2$, the simplex \eqref{s4 e1} is $\Delta_1 = [1,\bar e]$. Recall that 
\[\Delta_1^{(1)}=(1,\bar{e}],\quad \Delta_1^{(2)}=[1,\bar{e}),  \quad F_1=\{1\},\quad F_2=\{\bar{e}\}.\]
The subsets $\E_1=\pi^{-1}(\Delta_1^{(1)})$ and $\E_2=\pi^{-1}(\Delta_1^{(2)})$ of $\E$ are given by 
\[\E_1=\{E\in \E: r_1(E)< r_2(E) \},\quad \E_2=\{E\in \E: r_2(E)<r_3(E)\}.\]
We have the Mayer-Vietoris sequence for the decomposition $\E=\E_1\cup\E_2$:
\beq\label{equ:2:mv}
\cdots\to H_k(\E_1)\oplus H_k( \E_2)\to H_k(\E)\to H_{k-1}(\E_1\cap \E_2)\to H_{k-1}(\E_1)\oplus H_{k-1}( \E_2) \to \cdots
\eeq
In order to prove Theorem \ref{thm:top} by using \eqref{equ:2:mv},
we compute the homology groups of $\E_1$, $\E_2$ and $\E_1\cap \E_2$.

The following lemma helps us to simplify the computation of the homology groups of $\E_1$, $\E_2$ and $\mc{W}_{12}=\E_1\cap \E_2$. 

\begin{lemma}\label{lem:2:retract} The following statements hold:
\begin{enumerate}
\item For $i=1,2$, $\pi^{-1}(F_{3-i})$ is a deformation retract of $\E_i$. 
\item For any given point $P \in \Delta_1^{(1)}\cap \Delta_1^{(2)}$, $\pi^{-1}(P)$ is a deformation retract of $\E_1\cap \E_2$.
\end{enumerate}
As a direct consequence, we have for all $k\geq 1$, 
$$
{\begin{split}
  H_k(\E_i) &=H_k(\pi^{-1}(F_{3-i})), \ \ \  \  i=1,2,\\
  H_k(\E_1\cap \E_2) &=H_k(\pi^{-1}(P)),\ \ \ P \in \Delta_1^{(1)}\cap \Delta_1^{(2)}.  
 \end{split}}
 $$
\end{lemma}

\begin{proof}
For part (1), it suffices to consider the case when $i=1$. Recall that $\Delta_1^{(1)}=(1,\bar{e}]$ and $F_2=\{\bar{e}\}$. 
Let $G: \Delta_1^{(1)}\times [0,1]\to \Delta_1^{(1)}$ be a  strong deformation retraction of $\Delta_1^{(1)}$ onto $F_2$ 
(such deformation clearly exists). 
We then define $\mc{G}:\E_1\times [0,1]\to \E_1$ as follows. 
For any $E\in \E_1$, let $\mc{G}(E,t)$ be the ellipsoid such that its axial directions are all the same with $E$, 
and its axial lengths  are determined by
\beq\label{lem:2:e1}
(r_1(t),r_{2}(t),r_3(t))= \pi_2^{-1}\circ G(\pi(E),t).
\eeq
Namely, we continuously deform the axial lengths of $E$ while keeping the directions of its axes so that the resulting ellipsoid belongs to $\pi^{-1}(F_2)$. 
It is easy to check that $\mathcal G$ is a deformation retraction from $\E_1$ onto $\pi^{-1}(F_2)$.

For part (2), the argument is similar. 
\end{proof}

Since $\pi^{-1}(F_1)=\{E\in \E:r_1(E)=r_2(E)\}$ and $\pi^{-1}(F_2)=\{E\in \E:r_2(E)=r_3(E)\}$,
we see that $\pi^{-1}(F_1)$ and $\pi^{-1}(F_2)$ are both homeomorphic to $RP^2$. Using Lemma \ref{lem:2:retract} and the homology of $RP^2$, we have 
\begin{equation}\label{equ:2:e2}
H_k(\E_i)=0,  \  \ \forall \ k\ge3.
\end{equation}

Next, we study the homology groups of $\mc{W}_{12}=\E_1\cap \E_2$.

\begin{lemma}\label{lem:2:top} 
For any given $\r=(r_1,r_2,r_3)\in\mathbb L_{3}$ with $r_1<r_2<r_3$, denote $\E_{\r} = \pi_1^{-1}(\r)$, i.e., 
\[\E_{\r}=\{E\in \E: r_i(E)=r_i, 1\leq i\leq 3\}.\]Then $\E_{\r}\simeq SO(3)/\Gamma$, 
where $\Gamma=  \big\{A\in SO(3): A=\text{diag}\{\pm1,\pm1,\pm 1\}\big\}$ 
is a finite discrete subgroup of $SO(3)$.  It follows that $\E_{\r}$ is orientable and so
\[
H_{3}(\E_{\r}) = \mathbb Z.
\]
\end{lemma} 

\begin{proof}
We consider the following Lie group action on $\E_{\r}$:
\[
SO(3)\times\E_{\r}\to\E_{\r}, \ (g,E)\mapsto T_g(E),
\]
where $T_g$ is the linear transformation of $\R^{3}$ associated to $g$.
This groups action is clearly transitive and  the stabiliser of $E_{\r}=\{x\in\R^3:\sum_{i=1}^3x_i^2/r_i^2\le 1\}$ is given by
$ \Gamma $.
The subgroup $\Gamma$ is isomorphic to the dihedral group $D_2=\mb{Z}_2\oplus \mb{Z}_2$.
It follows that $\E_{\r}$ is  homeomorphic to $SO(3)/\Gamma$ \cite{Lee}. 
As $\Gamma$ is finite, we find that $SO(3)/\Gamma$ is an orientable manifold.
This together with $\dim(SO(3)/\Gamma)=3$ yields
$H_3(\E_{\r})=H_3(SO(3)/\Gamma)=\mb{Z}.$
\end{proof}

For any $P\in \Delta_1^{(1)}\cap \Delta_1^{(2)}$, 
let $\pi_2^{-1}(P)=(r_1,r_2,r_3)\in \mb{L}_3$ such that $r_1<r_2<r_3$. 
By Lemma \ref{lem:2:retract} and Lemma \ref{lem:2:top}, we obtain that 
\begin{equation}\label{equ:2:e12}
H_3(\E_1\cap \E_2)=\mb{Z}.
\end{equation}

Now, we are ready to prove Theorem \ref{thm:top} for $n=2$.

Taking $k=4$ in \eqref{equ:2:mv} and using \eqref{equ:2:e2}, we obtain the following short exact sequence
\beqs
0\to H_4(\E)\to H_3(\E_1\cap \E_2) \to0.
\eeqs
As a result, 
\[
H_4(\E) =  H_3(\E_1\cap \E_2).\]
This together with \eqref{equ:2:e12} gives $H_4(\E)=\mb{Z}$, and thus proves 
Theorem \ref{thm:top} for $n = 2$.

\subsection{Dimension $n=3$}\label{subsec:n3} 
{~}

When $n=3$, the simplex $\Delta_2=\{(x_2,x_3)\in \mb{R}^2: 1\leq x_2\leq x_3\leq \bar{e}\}$ is a triangle, and $F_i$, $1\leq i\leq 3$, are the sides of this triangle $\Delta_2$.  Recall that $\E_i=\pi^{-1}(\Delta_2^{(i)})$, $1\leq i\leq 3$, is given by 
\[\E_i=\{E\in \E: r_i(E)<r_{i+1}(E)\}.\]

Arguing similarly as in Lemma \ref{lem:2:retract}, we have the following result.
\begin{lemma}\label{lem:3:retract}
The following statements hold:
\begin{enumerate}
\item For $i=1,2,3$, $\pi^{-1}(P)$ is a deformation retract of $\E_i$, where $P=F_{j_1}\cap F_{j_2}$ and $j_1$ and $j_2$ are distinct such that $\{j_1,j_2\}=\{1,2,3\}\setminus \{i\}$.
\item For any point $P\in \bigcap_{k=1}^3 \Delta_2^{(k)}$, $\pi^{-1}(P)$ is a deformation retract of $\mc{W}_{123}=\E_1\cap \E_2\cap \E_3$.
\end{enumerate}
\end{lemma}

Suppose that $P= F_1\cap F_2$ and $\pi_2^{-1}(P)=(r_1,r_2,r_3,r_4)\in \mb{L}_4$. We then have $r_1=r_2=r_3$.  Hence, $\pi^{-1}(P)$ is homeomorphic to $RP^3$. By Lemma \ref{lem:3:retract}, we conclude that $\E_3$ is homotopy equivalent to $RP^3$. Similarly, we obtain that $\E_1$ is homotopy equivalent to $RP^3$ and $\E_2$ is homotopy equivalent to the Grassmannian $G(2,4)$. Using the homology of $RP^3$ and $G(2,4)$, we have for $1\leq i\leq 3$
\begin{equation}\label{equ:3:e1}
H_k(\E_i)=0, \quad \text{if } k\geq 5.
\end{equation}

By part (2) of Lemma \ref{lem:3:retract} and an analog of Lemma \ref{lem:2:top} for $n=3$ 
(see Lemma \ref{lem:top:mid} for the general case), 
we conclude that $\W_{123}$ is homotopy equivalent to $SO(4)/\Gamma$, 
where $\Gamma =  \big\{A\in SO(4): A=\text{diag}\{\pm1, \pm1, \pm1, \pm 1\}\big\}$ is a finite subgroup of $SO(4)$. 
Since $SO(4)/\Gamma$ is orientable and has dimension $6$, we have 
\begin{equation}\label{equ:3:e123}
H_6(\W_{123}) = H_6(SO(4)/\Gamma)=\mb{Z}.
\end{equation}

Next, we consider the homology groups of $\mc{W}_{12}=\E_1\cap \E_2$, $\mc{W}_{13}=\E_1\cap \E_3$ and $\mc{W}_{23}=\E_2\cap \E_3$. Recall that 
\[\mc{W}_{12}=\{E\in \E: r_1(E)<r_2(E)<r_3(E)\}.\]  
Take any point $P\in \Delta_2^{(1)}\cap \Delta_2^{(2)}\cap F_3$. 
Clearly, $\pi^{-1}(P)\in \mc{W}_{12}$. 
Arguing as in Lemma \ref{lem:2:retract}, 
we see that $\pi^{-1}(P)$ is a deformation retract of $\mc{W}_{12}$. 
Assume that $\pi_2^{-1}(P)=(r_1,r_2,r_3,r_4):={\r}$. 
We then have $r_1<r_2<r_3=r_4$. 
Denote $\E_{{\r}}=\pi^{-1}(P)$. 
We consider the Lie group action on $\E_{{\r}}$ as in Lemma \ref{lem:2:top}:  $SO(4)\times \E_{{\r}}\to \E_{{\r}}$. 
The  stabiliser of this group action is given by $S(O(1)\times O(1)\times O(2))$, 
i.e., the set of matrices in $O(1)\times O(1)\times O(2)$ with determinant 1. 
Therefore, we conclude that 
\[H_k(\mc{W}_{12})=H_k\big(SO(4)/S(O(1)\times O(1)\times O(2))\big),\, \forall\, k.\]
As the space $SO(4)/S(O(1)\times O(1)\times O(2))$ has dimension $5$, we get
\[H_k(\mc{W}_{12})=0,\,\forall\,k\geq 6.\]
The above discussion yields the following result.
\begin{lemma}\label{lem:3:retract2}
Let $1\leq j_1<j_2\leq 3$. Then 
\begin{enumerate}
\item $\pi^{-1}(P)$ is a deformation retract of $\mc{W}_{j_1j_2}$, where $P\in  \Delta_2^{(j_1)}\cap \Delta_2^{(j_2)}\cap F_i$
with $\{i\}=\{1,2,3\}\setminus \{j_1,j_2\}$.
\item $H_k(\mc{W}_{j_1j_2})=0$ for all $k\geq 6$.
\end{enumerate}
\end{lemma}

The Mayer-Vietoris sequence for the decomposition $\E=\E_1\cup(\E_2\cup \E_3)$ gives
\beq\label{equ:3:e3}
\cdots\to H_k(\E_1)\oplus H_k(\E_2\cup \E_3) \to H_k(\E)
\to H_{k-1}(\W_{12}\cup\W_{13} )\to H_{k-1}(\E_1)\oplus H_{k-1}(\E_2\cup \E_3)\cdots.
\eeq
For the decomposition $\E_{23}=\E_2\cup \E_3$, we have 
\beq \label{equ:3:e4}
\cdots\to 
H_k(\E_2)\oplus H_k(\E_3) \to H_k(\E_2\cup \E_3)
\to H_{k-1}( \W_{23})\to H_{k-1}(\E_2)\oplus H_{k-1}(\E_3) \to\cdots.
\eeq
Taking $k\ge7$ in \eqref{equ:3:e4} and using \eqref{equ:3:e1} and the part (2) in Lemma \ref{lem:3:retract2}, we see that 
\beq \label{equ:3:e5}
 H_k(\E_2\cup \E_3) = 0, \ \ \forall \ k\ge7.
\eeq
Letting $k=8$ in \eqref{equ:3:e3} and inserting \eqref{equ:3:e1} and \eqref{equ:3:e5} in the exact sequence, we obtain that 
\beq\label{equ:3:e6}
H_8 (\E)=  H_7( \W_{12}\cup\W_{13} ).
\eeq
The proof reduces to the computation of the right hand side of \eqref{equ:3:e6}, i.e., $H_7(\mc{W}_{12}\cup\mc{W}_{13})$.
The Mayer-Vietoris sequence 
\beqs
\begin{split}
\cdots \to H_7(\mc{W}_{12})\oplus H_7(\mc{W}_{13}) \to H_7( \mc{W}_{12}\cup \mc{W}_{13} ) 
\to H_6(\mc{W}_{123}) \to H_6(\mc{W}_{12})\oplus H_6(\mc{W}_{13})\to\cdots
\end{split}
\eeqs
together with the part (2) in Lemma \ref{lem:3:retract2} yields that
\begin{equation}\label{equ:3:e8}
H_7( \mc{W}_{12}\cup \mc{W}_{13} ) = H_6(\mc{W}_{123}).
\end{equation}

Combining \eqref{equ:3:e6}, \eqref{equ:3:e8} and \eqref{equ:3:e123}, we complete the proof of Theorem \ref{thm:top} for $n=3$.

\subsection{General dimensions}\label{subsec:n}
{~}

Now, we consider the general dimensions.
Before we use the Mayer-Vietoris sequence and the induction arguments, 
we first prove several lemmas concerning the homology groups of $\E_i$, $\E_{j_1j_2\cdots j_l}$ and $\mc{W}_{j_1j_2\cdots j_l}$
(see notations in \eqref{s4 e2} and \eqref{s4 e3}). 
Recall the following notations:
\begin{itemize} 
\item[(a)] If $P\in\Int\Delta_{n-1}$, then
\beqs
\pi^{-1}(P)\in \W_{12\cdots n} =\big \{E\in\E: r_i(E)\ne r_j(E) \ \text{whenever} \ i\ne j\big\}.
\eeqs
\item[(b)] If $P\in\big(\cap_{s=1}^l\Delta_{n-1}^{(j_s)} \big) \cap \big( \cap_{i\ne j_1, j_2,\cdots, j_l}F_i \big )$,
where $1\le l<n$ and $1\le j_1<\cdots<j_l\le n$,  then
\beqs
\pi^{-1}(P)\in\W_{j_1j_2\cdots j_l}\cap \{E\in\E: r_i(E)=r_{i+1}(E) \ \text{for all}\ i\ne j_s, \ s=1,\cdots,l\}.
\eeqs
\end{itemize}

The following lemma shows that $\pi^{-1}(P)$ in cases (a) and (b) above are deformation retracts of $\W_{12\cdots n}$
and $\W_{j_1j_2\cdots j_l}$, respectively.
It is the generalisation of Lemmas \ref{lem:2:retract}, \ref{lem:3:retract} and \ref{lem:3:retract2} for high dimensions. 

\begin{lemma}\label{lem:deform retract}
The two statements below hold.
\begin{itemize}
\item[(i)] For any given $P\in\Int\Delta_{n-1}$,
$\pi^{-1}(P)$ is a deformation retract of $\W_{12\cdots n}$.
\item[(ii)] For $1\le l<n$ and $1\le j_1<\cdots<j_l\le n$, if
$P\in\big(\cap_{s=1}^l\Delta_{n-1}^{(j_s)} \big) \cap \big( \cap_{i\ne j_1, j_2,\cdots, j_l}F_i \big )$.
Then $\pi^{-1}(P)$ is a deformation retract of $\W_{j_1j_2\cdots j_l}$.
\end{itemize}
\end{lemma}

\begin{proof}

For (i), let 
\[G:\ \Int\Delta_{n-1}\times[0,1]\to\Int\Delta_{n-1}\]
be a deformation retraction of $\Int\Delta_{n-1}$ onto $P$.
Define $\mathcal G:\W_{12\cdots n}\times[0,1]\to\W_{12\cdots n}$ as follows.
For any $E\in \W_{12\cdots n}$,
let $\mathcal G (E,t)$ be the ellipsoid such that its axial directions are all the same as $E$, 
and its axis lengths $r_i(t)$ are determined by
\beq\label{lem:deform retract e1}
(r_1(t),\cdots, r_{n+1}(t))= \pi_2^{-1}\circ G(\pi(E),t).
\eeq
It can be verified that $\mathcal G$ is a deformation retraction from $\W_{12\cdots n}$ onto $\pi^{-1}(P)$ that we want.

For (ii), the argument is similar.
Denote $W = \cap_{s=1}^l\Delta_{n-1}^{(j_s)}$.
Now let $G: W\times[0,1] \to W$ be a deformation retraction form $W$ onto $P$ such that
\begin{itemize}
\item $G(W\cap(\cap_{i\ne j_1,\cdots,j_l} F_i),t) \subset  W\cap(\cap_{i\ne j_1,\cdots,j_l} F_i)$ for all $t\in[0,1]$;
\item $G(Q,t) \in  W\setminus\cap_{i\ne j_1,\cdots,j_l} F_i $ for all $t\in[0,1)$
and $Q\in W\setminus\cap_{i\ne j_1,\cdots,j_l} F_i $.
\end{itemize}
We then define the deformation retraction $\mathcal G:\W_{j_1j_2\cdots j_l}\times [0,1]\to \W_{j_1j_2\cdots j_l}$ as follows:
$\mathcal G(E,t)$ keeps all the axis-directions of $E$ unchanged but
its axis-lengths $r_i(t)$ are again given by \eqref{lem:deform retract e1}.
\end{proof}

The next lemma gives the general case of Lemma \ref{lem:2:top}.
\begin{lemma}\label{lem:top:mid}
Suppose that $r_i$, $i=1, 2, \cdots, n+1$, 
are distinct positive constants such that $\r=(r_1,\cdots, r_{n+1})\in\mathbb L_{n+1}$.
Denote $\E_{\r} = \pi_1^{-1}(\r)$. Then $\E_{\r}$ is homeomorphic to $SO(n+1)/\Gamma$, 
where $\Gamma=  \big\{A\in SO(n+1): A=\text{diag}\{\pm1,\cdots,\pm 1\}\big\}$ is a finite discrete subgroup of $SO(n+1)$.
As a result $\E_{\r}$ is orientable and has dimension $n^*$, and 
\beqs
H_{n^*}(\E_{\r}) = \mathbb Z.
\eeqs
\end{lemma}

\begin{proof}
Consider the Lie group action on $\E_{\r}$:
\beq\label{lem:top mid e2}
SO(n+1)\times\E_{\r}\to\E_{\r}, \ (g,E)\mapsto T_g(E),
\eeq
where $T_g$ represents the linear transformation of $\R^{n+1}$ associated to $g$.
Then the stabiliser of
\beq\label{lem:top mid e1}
E_{\r}=\Big\{x\in\R^{n+1}:\sum_{i=1}^{n+1}x_i^2/r_i^2\le 1\Big\}
\eeq
is given by
$
\Gamma =  \big\{A\in SO(n+1): A=\text{diag}\{\pm1,\cdots,\pm 1\}\big\},
$
the set of diagonal $(n+1)\times(n+1)$ matrices with determinant $1$. 
Since the group action is transitive and $\Gamma$ is closed in $SO(n+1)$,
we conclude that $\E_{\r}$ has a smooth manifold structure such that
\beqs
\mathcal F:SO(n+1)/\Gamma\to \E_{\r}, \ \mathcal F(g\Gamma)\mapsto T_g(E_{\r}),
\eeqs
is a diffeomorphism \cite{Lee}.
 As such manifold structure yields the same topology of $\E_{\r}$ induced by the Hausdorff metric,
 we see that $\E_{\r}\simeq SO(n+1)/\Gamma$ and these two spaces have the same homology.

Since $\Gamma$ is a finite discrete group and the orientation of $SO(n+1)$
is preserved by all the diffeomorphisms of $\Gamma$, 
we conclude that $SO(n+1)/\Gamma$ is an orientable closed manifold.
As $\dim(SO(n+1)/\Gamma) = n^*$, we see that (\cite{Hat})
\beqs
H_{n^*}(SO(n+1)/\Gamma) = \mathbb Z.
\eeqs
This completes the proof.
Indeed, one can also show that $SO(n+1)/\Gamma$ is diffeomorphic to the complete flag variety in $\mb{R}^{n+1}$
(see \cite{Br05}).
\end{proof}

The following result is corresponding to  the general case of Lemma \ref{lem:3:retract2}.

\begin{lemma}\label{lem:top:facet}
Suppose that $\r=(r_1,\cdots,r_{n+1})\in\mathbb L_{n+1}$ such that 
$r_{i_k+1}$ has multiplicity $m_k$ for $k=1, \cdots, l$ and all other components  ${r_i}$ are single,
where 
$${\begin{split}
0\le i_1<i_2 & <\cdots<i_l\le n-1, \ \ i_j+m_j\le i_{j+1},\ \   i_l+m_l\le n+1,\ \ 
  {\Small\text{${\sum}_{j=1}^l $}} m_j\le n+1.
 \end{split}} $$  
 Let $\E_{\r} = \pi_1^{-1}(\r)$. Then 
\beq\label{lem:top facet e0}
H_k(\E_{\r}) = 0, \ \ \text{if} \ k\ge n^*+1- {\Small\text{${\sum}_{j=1}^l  \frac{m_j(m_j-1)}{2} $}}.
\eeq
\end{lemma}
\begin{proof} 
The components of $\r$ can be divided into two groups:
\[\{r_{s_1}, \cdots, r_{s_p}\}\quad \text{and}\quad \{r_{i_1+1},\cdots, r_{i_l+m_l}\}.\]
Components in the first group are single ones  and components in the second are multiple ones (see Remark \ref{r4.1}).
We have $p+\sum_{j=1}^l m_j=n+1$.

For simplicity,  we may assume that (after a proper permutation) $\r$  can be written as
\beq \label{vecr}
(r_{s_1}, \cdots, r_{s_p}\ |\ r_{i_1+1},\cdots r_{i_1+m_1},\cdots, r_{i_l+1}, \cdots, r_{i_l+m_l}) , 
\eeq
where the components before the symbol $|$ are single ones and the components after the symbol $|$ are multiple ones,
and the components are in the ascending order
$r_{s_1}<\cdots < r_{s_p}$, $r_{i_1+1}<r_{i_2+1}<\cdots<r_{i_l+1}$.
  
Consider the Lie group action on $\E_{\r}$ as in \eqref{lem:top mid e2}.
This group action is transitive and the stabiliser of $E_{\r}$ (given by \eqref{lem:top mid e1}) 
is the collection of matrices in the form
\beqs
\mathcal S=\text{diag}\, \{\pm 1, \cdots, \pm1\ |\ \mathcal O_{m_1},\cdots, \mathcal O_{m_l}\}
\eeqs
with the property $\det \mathcal S=1$,
where $\mathcal O_{m_k}\in O(m_k)$, the set of $m_k\times m_k$ orthogonal matrices.
By the same argument as in Lemma \ref{lem:top:mid}, we obtain
\beqs
 \E_{\r}\simeq SO(n+1)/S\big(O(1)\times\cdots \times O(1)\times O(m_1)\times  \cdots \times O(m_l)\big).
\eeqs
It is known that the space on the right-hand side has dimension $n^*-\sum_{j=1}^lm_j(m_j-1)/2$.
We then deduce \eqref{lem:top facet e0} as desired. 
\end{proof}


By using Lemma \ref{lem:deform retract} and Lemma \ref{lem:top:mid}, we have the following conclusion. 
\begin{lemma}\label{lem:n:Y}
We have $H_{n^*}(\mathcal W_{12\cdots n}) = \mathbb Z$.
\end{lemma}
\begin{proof}
Let $P$ be a point of $\cap_{i=1}^n\Delta_{n-1}^{(i)}$. By Lemma \ref{lem:deform retract}
\beqs
H_*(\mathcal W_{12\cdots n}) = H_*(\pi^{-1}(P)).
\eeqs
Since $\pi^{-1}_2(P) = \r\in \mathbb L_{n+1}$ satisfies $r_1<r_2<\cdots<r_{n+1}$, 
it follows from Lemma \ref{lem:top:mid} that
\beqs
H_{n^*}(\pi^{-1}(P)) = \mathbb Z.
\eeqs
This completes the proof.
\end{proof}

\begin{lemma}\label{lem:n:Y0}
 For any $1\leq j_1<j_2\cdots <j_l\leq n$, we have 
\[H_k(\mathcal W_{j_1j_2\cdots j_l})=0, \quad \forall\, k\ge n^*+1-  {\Small\text{${\sum}_{s=1}^{l+1} \frac{m_s(m_s-1)}{2} $}},\]
where $m_s = j_s-j_{s-1}$ and $j_0=0$, $j_{l+1} =n+1$.
\end{lemma}
\begin{proof}
Let $P$ be a point in $(\cap_{s=1}^l\Delta_{n-1}^{(j_s)})\cap (\cap_{i\ne j_1, j_2,\cdots, j_l}F_i)$.
By Lemma \ref{lem:deform retract}, we find
\beq\label{lem:n:Y0:e1}
 H_*(\mathcal W_{j_1j_2\cdots j_l}) = H_*(\pi^{-1}(P)).
\eeq
As $\pi^{-1}_2(P)=\r\in \mathbb L_{n+1}$ satisfies
\beqs
& r_1=\cdots = r_{j_1}<r_{j_1+1}=\cdots = r_{j_2}<r_{j_2+1} =\\
& \cdots
  =r_{j_{l-1}}<r_{j_{l-1}+1} \cdots=r_{j_l}<r_{j_l+1} = \cdots = r_{n+1},
\eeqs
we obtain the conclusion by Lemma \ref{lem:top:facet} and \eqref{lem:n:Y0:e1}. 
\end{proof}

Propositions \ref{prop:key0} and \ref{prop:key1} below are consequences of Lemmas \ref{lem:n:Y} and \ref{lem:n:Y0},
which can be viewed as a generalisation of these two lemmas.

\begin{proposition}\label{prop:key0}
Suppose $1\le p_1<\cdots <p_r<j\le n$, we have 
\beq\label{prop:key0:e0}
H_{k}(\cup_{l=j}^n \mathcal W_{p_1\cdots p_r l})=
 0 , \ \ \ \text{if} \ k \ge n^*+n+1-j-  {\Small\text{${\sum}_{s=1}^{r+1}\frac{m_s(m_s-1)}{2} $}} ,
\eeq
where $m_s = p_s-p_{s-1}$, $p_0 = 0$ and $p_{r+1} =j $.
Furthermore, if $k\ge n^*+n+1-j-\frac{j(j-1)}{2}$, then
\beq\label{prop:key0:e1}
H_k(\E_{j(j+1)\cdots n})=0.
\eeq
\end{proposition}

\begin{proof}
For $j=n$, \eqref{prop:key0:e0} follows from Lemma \ref{lem:n:Y0}.


We now verify \eqref{prop:key0:e0} by the induction argument on $j$.
For this purpose, let us assume that \eqref{prop:key0:e0} holds when $j=n-m$ for some $m\ge0$.
We next show that \eqref{prop:key0:e0} holds for $j=n-(m+1)$.
The Mayer-Vietoris sequence for the decomposition
\beqs
 \cup_{l=n-m-1}^{n}\W_{p_1\cdots p_rl } 
 =\W_{p_1\cdots p_r(n-m-1) } \cup (\cup_{l=n-m}^{n}\W_{p_1\cdots p_rl } )
\eeqs
yields
\beq\label{prop:key0:e2}
\begin{split}
\cdots & \to H_k(\cup_{l=n-m}^{n}\W_{p_1\cdots p_rl } )\oplus H_k(\W_{p_1\cdots p_r(n-m-1) })   \\
 & \to H_k( \cup_{l=n-m-1}^{n}\W_{p_1\cdots p_rl } ) \to H_{k-1}(\cup_{l=n-m}^{n}\W_{p_1\cdots p_r (n-m-1) l } ) \\
 & \to H_{k-1}(\cup_{l=n-m}^{n}\W_{p_1\cdots p_rl } )\oplus H_{k-1}(\W_{p_1\cdots p_r(n-m-1) })\to \cdots.
\end{split}
\eeq

By our induction assumption, \eqref{prop:key0:e0} holds when $j=n-m$. That is
\beq\label{prop:key0:e3}
 H_k(\cup_{l=n-m}^{n}\W_{p_1\cdots p_rl } )=0,  \ \  \text{if} \ k\ge k(m),
\eeq
where $k(m)$ is an integer function of $m$ given by 
\[k(m):=n^*+m+1-  {\Small\text{${\sum}_{s=1}^r\frac{m_s(m_s-1)}{2} $}} -  {\Small\text{$ \frac{(n-m-p_r)(n-m-p_r-1)}{2} $}}. \]
By Lemma \ref{lem:n:Y0},
\beq\label{prop:key0:e4}
H_k(\W_{p_1\cdots p_r(n-m-1) }) =0, \ \ \text{if} \ k\ge k'(m),
\eeq
where $k'(m)$ is another integer function of $m$ given by 
\[
k'(m):=n^*+1-  {\Small\text{${\sum}_{s=1}^r\frac{m_s(m_s-1)}{2} $}}
                  -  {\Small\text{$ \frac{(n-m-p_r-1)(n-m-p_r-2)}{2}-\frac{(m+2)(m+1)}{2} $}}.\]

It can be verified that $k(m+1)\ge\max\{k(m),k'(m)\}+1$. 
Now inserting \eqref{prop:key0:e3} and \eqref{prop:key0:e4} in the long exact sequence \eqref{prop:key0:e2}, we obtain then
\beq\label{prop:key0:e5}
H_k( \cup_{l=n-m-1}^{n}\W_{p_1\cdots p_rl } )=  H_{k-1}(\cup_{l=n-m}^{n}\W_{p_1\cdots p_r (n-m-1) l } ), \ \ \text{if} \ k\ge k(m+1).
\eeq
By our induction assumption again, the right hand side above
\beqs
H_{k-1}(\cup_{l=n-m}^{n}\W_{p_1\cdots p_r (n-m-1) l } )=0, \ \ \text{if} \ k\ge k(m+1).
\eeqs
Hence, by \eqref{prop:key0:e5}, we conclude that \eqref{prop:key0:e0} holds when $j=n-(m+1)$.

Note that $\W_i = \E_i$ and $\E_{j(j+1)\cdots n} = \cup_{l=j}^n\W_l$.
By the same discussion as above but deleting $p_i$'s, we obtain \eqref{prop:key0:e1}.
\end{proof}

\begin{proposition}\label{prop:key1}
For any $2\leq j\leq n$, we have 
\beq\label{prop:key1:e0}
H_{n^*+n-j}(\cup_{l=j}^n \W_{12\cdots (j-1)l})=\mb{Z}.
\eeq
In particular, if $j=2$, then $H_{n^*+n-2}( \cup_{l=2}^n \mathcal W_{1l} )=\mb{Z}$.
\end{proposition}

\begin{proof}
For $j=n$, \eqref{prop:key1:e0} is the conclusion of Lemma \ref{lem:n:Y}.
Suppose by induction argument that \eqref{prop:key1:e0} holds for $j=n-m$ for some $m\ge0$.
Applying the Mayer-Vietoris sequence to the pair 
$\cup_{l=n-m}^{n}\W_{12\cdots (n-m-2) l } $ and $\W_{12\cdots (n-m-2)( n-m-1) } $,
we obtain
\beqs
\begin{split}
\cdots & \to H_k(\cup_{l=n-m}^{n}\W_{12\cdots (n-m-2) l } )\oplus H_k(\W_{12\cdots (n-m-2)(n-m-1)})\\
  & \to H_k(\cup_{l=n-m-1}^{n}\W_{12\cdots (n-m-2) l }  )\to 
          H_{k-1}(\cup_{l=n-m}^{n}\W_{12\cdots (n-m-1) l }  ) \\
   &\to H_{k-1}(\cup_{l=n-m}^{n}\W_{12\cdots (n-m-2) l } )\oplus H_{k-1}(\W_{12\cdots (n-m-2)(n-m-1)})\to \cdots.
\end{split}
\eeqs
It follows from \eqref{prop:key0:e0} in Proposition \ref{prop:key0} and Lemma \ref{lem:n:Y0} that
\beqs
 H_k(\cup_{l=n-m}^{n}\W_{12\cdots (n-m-2) l } )= H_k(\W_{12\cdots (n-m-2)(n-m-1)})=0, \ \ \text{if} \ k \ge n^*+m,
\eeqs
and therefore the long exact sequence above implies that
\beqs
H_{n^*+m+1}(\cup_{l=n-m-1}^{n}\W_{12\cdots (n-m-2) l }  ) = H_{n^*+m}(\cup_{l=n-m}^{n}\W_{12\cdots (n-m-1) l }).
\eeqs
Hence, \eqref{prop:key1:e0} follows when $j=n-m-1$ by our induction assumption.
This completes the proof.
\end{proof}

Now, we are ready to give the proof of Theorem \ref{thm:top} for general dimensions. 
\begin{proof}[Proof of Theorem \ref{thm:top}]
The Mayer-Vietoris sequence for the decompostion $\E=\E_1\cup \E_{23\cdots n}$ implies
\beq\label{thm:MV}
\begin{split}
\cdots & \to H_{n^*+n-1}(\E_1)\oplus H_{n^*+n-1}(\E_{23\cdots n}) \to H_{n^*+n-1}(\E)\\
 & \to H_{n^*+n-2}(\E_1\cap \E_{23\cdots n})
\to  H_{n^*+n-2}(\E_1)\oplus H_{n^*+n-2}(\E_{23\cdots n}) \to \cdots.
\end{split}
\eeq
By virtue of \eqref{prop:key0:e1} in Proposition \ref{prop:key0} and Lemma \ref{lem:n:Y0} (for $l=1$ and $j_1=1$),
\beqs
H_{n^*+n-1}(\E_1)=H_{n^*+n-1}(\E_{23\cdots n})= H_{n^*+n-2}(\E_1)= H_{n^*+n-2}(\E_{23\cdots n})=0.
\eeqs
Hence, by \eqref{thm:MV},
\beqs
H_{n^*+n-1}(\E)  = H_{n^*+n-2}(\E_1\cap \E_{23\cdots n}).
\eeqs
Since $\E_1\cap \E_{23\cdots n}= \cup_{l=2}^n\W_{1l}$, we complete the proof by Proposition \ref{prop:key1}.
\end{proof}

\begin{remark}\label{r4.12}
For any given $ k\ge n^*+n$, 
by Proposition \ref{prop:key0} and Lemma \ref{lem:n:Y0} we have
\beqs
H_{i}(\E_1) = H_{i}(\E_{23\cdots n})  = 0\ \ \text{for} \ i= k-1\ \text{or}\ k.
\eeqs
Using \eqref{thm:MV} with $n^*+n-1$ replaced by $k$, we then obtain
\beqs
 H_k(\E) = H_{k-1}(\E_1\cap \E_{23\cdots n}) = H_{k-1}(\cup_l^n\W_{1l}).
\eeqs
By \eqref{prop:key0:e0} (with $j=2$), the right hand side above $H_{k-1}(\cup_l^n\W_{1l})=0$. Therefore
\beq\label{large k hom}
 H_k(\E) =0 \ \text{for all} \ k\ge n^*+n.
\eeq
\end{remark}


\section{Proof of Theorem \ref{thm:estimates}}\label{sec:ap}

In this section, we prove Theorem \ref{thm:estimates} by showing 
\begin{itemize}
\item [(i)] the Gauss curvature of $\M_t$ is bounded from above,
\item [(ii)] the principal curvatures of $\M_t$ have a positive lower bound.
\end{itemize}
By approximation, we may assume directly that $f$ is $C^2$-smooth.
The Gauss curvature flow has been extensively studied.
The technique and calculation presented here are similar to those in \cite{LSW20}.

Let $X(\cdot,t)$ be the solution of the flow \eqref{flow}.
Recall that the Gauss curvature of $X(\cdot,t)$ is given by \eqref{s2 e1},
and the principal radii of curvature of $X(\cdot,t)$ are eigenvalues of the matrix $\{b_{ij}\}$, 
where
\beqs
b_{ij} = u_{ij}+u\delta_{ij},
\eeqs
where $u$ is the support function of $X(\cdot,t)$.

First, we derive an upper bound for the Gauss curvature. 

\begin{lemma}\label{lem:ap:K}
Let $X(\cdot, t)$ be a uniformly convex solution to the flow \eqref{flow} for $t\in [0,T)$. 
Suppose that the support function $u$ satisfies \eqref{thm:estimates:C0}.
Then there exists a constant $ C$ depending on $n, p, \min_{\S^n}f$, $\max_{\S^n} f$, 
the initial condition $\M_0$, and the constant $C_0$ in   \eqref{thm:estimates:C0}, such that
\beq\label{KC}
K(\cdot, t)\leq C,\quad \forall\, t\in [0,T).
\eeq
\end{lemma}

\begin{proof} 
We introduce the auxiliary function
\[Q=-\frac{u_t}{u-\eps_0} = \frac{Ku^{p}f-u}{u-\eps_0},\]
where $\eps_0=\frac{1}{2}\min_{\mb{S}^n\times [0,T)}u>0$.
It suffices to show that $Q(x,t)\le C$ $\forall$ $(x,t)\in\S^n\times[0,T)$.

For any given $T'\in (0, T)$, we
assume that $Q$ attains its maximum over $\S^n\times[0,T']$ at $(x_0,t_0)$.
If $t_0 = 0$, then $\max_{\S^n\times[0,T']}Q= \max_{\S^n}Q(\cdot,0)$ and we are through.
If $t_0>0$, then at the point $(x_0,t_0)$, we have
\begin{equation}\label{equ:ap:K:1}
0=\nabla_i Q=-\frac{u_{ti}}{u-\eps_0}+\frac{u_tu_i}{(u-\eps_0)^2}.
\end{equation}
Hence $u_{ti} = -Qu_i $ and we have
\beqn\label{equ:ap:K:2}
0\geq \nabla_{ij}^2 Q& =&
-\frac{u_{tij}}{u-\eps_0}+\frac{u_{ti}u_j+u_{tj}u_i+u_tu_{ij}}{(u-\eps_0)^2}-\frac{2u_t u_iu_j}{(u-\eps_0)^3} \notag \\ 
& =&-\frac{u_{tij}}{u-\eps_0} +\frac{u_tu_{ij}}{(u-\eps_0)^2}.
\eeqn
It follows that 
\begin{equation}\label{equ:ap:K:b}
-b_{ijt}=-u_{ijt}-u_t \delta_{ij}\leq (b_{ij}-\eps_0\delta_{ij})Q.
\end{equation}
Let $\{h^{ij}\}$ be the inverse matrix of $\{b_{ij}\}$. Then
\beqs
  {\Small\text{$ \sum \, $}} h^{ii} \ge n\big(  {\Small\text{$ \prod $}} h^{ii}\big)^{\frac1n} = n K^{1/n}.
\eeqs
This, together with \eqref{equ:ap:K:b},  yields
\beq\label{equ:ap:K:4}
{\begin{split}
\p_t K 
    & = -K  {\Small\text{$ \sum \, $}} h^{ij} b_{ijt} \\
    & \le (n-\eps_0  {\Small\text{$ \sum \, $}}h^{ii})KQ\\
    & \le CQ^2-\frac{\eps_0}{C}Q^{2+1/n}.
    \end{split}}
\eeq
We next compute, at $(x_0,t_0)$,
\beqn\label{equ:ap:K}
0\le \p_t Q & =&-\frac{u_{tt}}{u-\eps_0}+Q^2 \notag \\
& =&\frac{1}{u-\eps_0}\frac{\p}{\p t}\Big( Ku^pf\Big) +Q+Q^2\notag \\ 
& \le&\frac{1}{u-\eps_0}\Big( fu^p\p_t K \Big)+CQ^2,
\eeqn
where we assume without loss of generality that $K\approx Q\gg 1$. 

Combining \eqref{equ:ap:K:4} and \eqref{equ:ap:K}, we obtain, at $(x_0,t_0)$,
\beqs
0\leq \big(C-\eps_0Q^{1/n} \big) Q^2.
\eeqs
This implies that $\max_{\S^n\times[0,T']} Q$ is bounded from above. 
As this bound is independent of $T'$, by sending $T'\to T$, we complete the proof. 
\end{proof}

Next, we derive a lower bound on the principal curvatures.

\begin{lemma}\label{lem:ap}
Let $X(\cdot, t)$ be a uniformly convex solution to the flow \eqref{flow} for $t\in [0,T)$. 
Assume the support function $u$ satisfies \eqref{thm:estimates:C0}.
Then there exists a constant $\bar{\kappa}$ depending on $n, p, C_0, \min_{\S^n}f, \|f\|_{C^{1,1}(\S^n)}$,
and the initial condition $\M_0$, such that 
\beq\label{kappai}
\kappa_i(\cdot, t)\geq \bar{\kappa} \quad \forall\, t\in [0,T),\, 1\leq i\leq n,
\eeq
where $\kappa_i$'s are the principal curvatures of $X(\cdot,t)$.
\end{lemma}
\begin{proof}
Consider the following auxiliary function
\[\wti{w}(x,t)=\log \lambda_{\max}(\{b_{ij}\})-A\log u+B|\nabla u|^2, \]
where $A$ and $B$ are large constants to be determined, 
and $\lambda_{\max}(\{b_{ij}\})$ denotes the maximal eigenvalue of $\{b_{ij}\}$.
Our purpose is to show that $\wti{w}$ is bounded from above.

For any given $T'\in (0, T)$,
assume that $ \wti{w}(x,t)$ achieves its maximum over $\mb{S}^n\times[0,T']$ at some point $(x_0, t_0)$.
We also suppose $t_0>0$,
otherwise estimate \eqref{kappai} follows from the initial condition.
By a proper rotation, 
we may assume that $\{b_{ij}\}$ is diagonal at $(x_0,t_0)$ and $\lambda_{\max}(\{b_{ij}\})(x_0,t_0)=b_{11}(x_0,t_0)$.

Then the function 
\[w(x,t)=\log b_{11}-A\log u+B|\nabla u|^2\]
attains its maximum at $(x_0,t_0)$.  
We may assume $b_{11}\gg 1$, otherwise we are through.
Denote by $\{h^{ij}\}$ the inverse matrix of $\{b_{ij}\}$. At $(x_0,t_0)$, we have 
\beqn\label{equ:pc:max:1}
0=\nabla_i w&=&h^{11}\nabla_i b_{11} -A\frac{u_i}{u}+2B  {\Small\text{$ \sum_k $}}\, u_ku_{ki} \notag\\
&=& h^{11}(u_{i11}+u_1\delta_{i1}) -A\frac{u_i}{u}+2B  u_iu_{ii},
\eeqn
and 
\begin{equation}\label{equ:pc:max:2}
0\geq \nabla_{ii}w=h^{11}\nabla_{ii}^2 b_{11}-(h^{11})^2(\nabla_i b_{11})^2-A\Big(\frac{u_{ii}}{u}-\frac{u_i^2}{u^2}\Big) +2B\Big(u_{ii}^2 + {\Small\text{$ \sum_k$}}  u_k u_{kii}\Big).
\end{equation}
In the above, we have used the properties that 
$\nabla_k b_{ij}$ are symmetric in all indices and that $\nabla_k b^{ij}=-h^{il}h^{jp}\nabla_k b_{lp}$.

We also have 
\[
\p_t w 
=h^{11} (u_{11t} +u_t) -A\frac{u_t}{u}+2B  {\Small\text{$ \sum $}}\, u_k u_{kt}.
\]
Next, we estimate the term $b^{11}u_{11t}$. Recall that 
\begin{equation}\label{equ:pc:u_t}
\log (u-u_t)=\log K +\log \big(fu^p\big).
\end{equation}
Set 
\[\phi(x,u)=\log \big( fu^p \big).\]
Differentiating \eqref{equ:pc:u_t} gives 
\beqn\label{lem:ap:ukt}
\frac{u_k-u_{kt}}{u-u_t}&=&-  {\Small\text{$ \sum$}}\, h^{ij}\nabla_k b_{ij}+\nabla_k \phi\notag \\
&=&-  {\Small\text{$ \sum \, $}}h^{ii} \big(u_{kii}+u_i \delta_{ik}\big) +\nabla_k \phi, 
\eeqn
and 
\begin{equation}\label{equ:pc:u_t11}
\frac{u_{11}-u_{11t}}{u-u_t}-\frac{(u_1-u_{1t})^2}{(u-u_t)^2}
=-  {\Small\text{$ \sum \, $}}h^{ii}\nabla_{11}^2 b_{ii} +  {\Small\text{$ \sum \, $}}h^{ii}h^{jj}(\nabla_1 b_{ij})^2 +\nabla_{11}^2\phi.
\end{equation}
By \eqref{equ:pc:u_t11} and the Ricci identity $\nabla_{11}^2 b_{ii}=\nabla_{ii}^2 b_{11}-b_{11}+b_{ii}$, we have 
\beqn\label{equ:pc:w}
\frac{\p_t w}{u-u_t}& =&
h^{11}\Big[\frac{u_{11t}-u_{11}}{u-u_t}+\frac{u_{11}+u-u+u_t}{u-u_t}\Big] 
-\frac{A}{u}\frac{u_t-u+u}{u-u_t}+2B\frac{\sum u_k u_{kt}}{u-u_t} \notag \\
&\le &h^{11}\Big[{\Small\text{$ \sum \, $}}h^{ii}\nabla_{11}^2 b_{ii} -{\Small\text{$ \sum \, $}}h^{ii}h^{jj}(\nabla_1 b_{ij})^2 -\nabla_{11}^2\phi \Big] \notag\\ 
&& +\frac{1}{u-u_t}+\frac{A}{u}-\frac{A}{u-u_t}+2B\frac{\sum u_k u_{kt}}{u-u_t} \notag \\ 
& \leq & h^{11}\Big[{\Small\text{$ \sum \, $}}h^{ii} \big(\nabla_{ii}^2  b_{11}-b_{11}+b_{ii}\big)-{\Small\text{$ \sum \, $}}h^{ii}h^{jj}(\nabla_1 b_{ij})^2\Big]
-h^{11}\nabla_{11}^2\phi  \\ 
&& +\frac{1-A}{u-u_t}+\frac{A}{u}+ 2B\frac{\sum u_k u_{kt}}{u-u_t}. \notag
\eeqn
Inserting \eqref{equ:pc:max:1} and \eqref{equ:pc:max:2} into \eqref{equ:pc:w}, we obtain, at $(x_0,t_0)$,
\beqs
\frac{\p_t w}{u-u_t}&\leq&
 {\Small\text{$ \sum \, $}}h^{ii}\Big[(h^{11})^2 (\nabla_i b_{11})^2 +A\Big(\frac{u_{ii}}{u}-\frac{u_i^2}{u^2}\Big)  
-2B\Big(u_{ii}^2 +{\Small\text{$ \sum \, $}}u_k u_{kii}\Big)\Big] \\ 
&& +h^{11}{\Small\text{$ \sum \, $}}h^{ii}(b_{ii}-b_{11}) -h^{11}{\Small\text{$ \sum \, $}}h^{ii}h^{jj}(\nabla_1 b_{ij})^2-h^{11}\nabla_{11}^2\phi  \\ 
&& +\frac{1-A}{u-u_t}+\frac{A}{u}+ 2B\frac{\sum u_k u_{kt}}{u-u_t}\\
&\le&{\Small\text{$ \sum \, $}}h^{ii}\Big[ A\Big(\frac{u_{ii}}{u}-\frac{u_i^2}{u^2}\Big)  
-2B\Big(u_{ii}^2 +{\Small\text{$ \sum \, $}}u_k u_{kii}\Big)\Big] \\ 
&&-h^{11}\nabla_{11}^2\phi + 2B\frac{\sum u_k u_{kt}}{u-u_t}+\frac{1-A}{u-u_t}+CA\\
&\le&-A{\Small\text{$ \sum \, $}}h^{ii} -2B {\Small\text{$ \sum \, $}}b_{ii}-2B {\Small\text{$ \sum \, $}}h^{ii} u_k u_{kii} \\ 
&& -h^{11}\nabla_{11}^2\phi  + 2B\frac{\sum u_k u_{kt}}{u-u_t}+\frac{1-A}{u-u_t}+CA+CB,
\eeqs
where $\sum h^{ii}h^{11} (\nabla_i b_{11})^2 \le \sum h^{ii}h^{jj}(\D_1b_{ij})$ is used in the second inequality.
By \eqref{lem:ap:ukt}, 
\beqn\label{lem:ap:e1}
\frac{\p_t w}{u-u_t}&\le&(2B|\D u|^2-A){\Small\text{$ \sum \, $}}h^{ii} -2B {\Small\text{$ \sum \, $}}b_{ii} -h^{11}\nabla_{11}^2\phi  \notag\\ 
&&-2B{\Small\text{$ \sum \, $}}u_k\D_k\phi +\frac{2B|\D u|^2+1-A}{u-u_t}+CA+CB  \notag \\
&\le&(2B|\D u|^2-A){\Small\text{$ \sum \, $}}h^{ii} -2B {\Small\text{$ \sum \, $}}b_{ii}+Cb_{11}   \notag \\ 
&& +\frac{2B|\D u|^2+1-A}{u-u_t}+CA+CB.
\eeqn
Choose $B$ large such that $B\sum b_{ii}\ge Cb_{11}$, 
and let $A=2B \max_{\mb{S}^n\times [0,T)} |\nabla u|^2+1$.
Since $\p_tw\ge 0$ at $(x_0,t_0)$,
we obtain by \eqref{lem:ap:e1} that
\beqs
0\le \frac{\p_t w}{u-u_t}\le-B {\Small\text{$ \sum \, $}}b_{ii} +CA+CB.
\eeqs
Hence, $\lambda_{\max}(\{b_{ij}\})(x_0,t_0)$ is bounded and so $\max_{\S^n\times[0,T']}\lambda_{\max}(b_{ij})\le C$.
Since this upper bound is independent of $T'$, we then let $T'\to T$ and finish the proof.
\end{proof}

\begin{acknowledgements} 
The first author was supported by ARC DE210100535. 
The second author was supported by NSFC12031017. 
The third author was supported by ARC DP200101084.
In the preparation of this paper,  
we consulted several people on the homology of the topological space of  ellipsoids. 
We would like to take this opportunity to express our gratitude to their help.
In particular, we would like to thank Shi Wang for several helpful discussions.
\end{acknowledgements}



\end{document}